\documentclass{amsart}
\usepackage{amsmath}
\usepackage{amssymb}
\usepackage{amsthm}

\DeclareMathOperator{\id}{id}

\DeclareMathOperator{\tr}{tr}

\DeclareMathOperator{\dvol}{dvol}

\DeclareMathOperator{\Ric}{Ric}

\DeclareMathOperator{\Rm}{Rm}

\DeclareMathOperator{\End}{End}
\DeclareMathOperator{\Hom}{Hom}

\newcommand{\og}{\overline{g}}

\newcommand{\lp}{\langle}
\newcommand{\rp}{\rangle}
\newcommand{\lv}{\lvert}
\newcommand{\rv}{\rvert}

\newcommand{\kso}{\mathfrak{so}}

\newcommand{\contr}{\lrcorner}
\newcommand{\semiplus}{
  \mbox{$
  \begin{picture}(12.7,8)(-.5,-1)
  \put(2,0.2){$+$}
  \put(6.2,2.8){\oval(8,8)[l]}
  \end{picture}$}}

\newcommand{\charconstant}{characteristic constant}

\newcommand{\mA}{\mathcal{A}}

\newcommand{\mC}{\mathcal{C}}

\newcommand{\mE}{\mathcal{E}}

\newcommand{\mR}{\mathcal{R}}

\newcommand{\mT}{\mathcal{T}}

\newcommand{\mV}{\mathcal{V}}

\newcommand{\bD}{\mathbb{D}}

\newcommand{\bN}{\mathbb{N}}

\newcommand{\bR}{\mathbb{R}}

\newcommand{\J}{\mathrm{J}}
\newcommand{\bg}{\mathbf{g}}

\newtheorem{thm}{Theorem}[section]
\newtheorem{prop}[thm]{Proposition}
\newtheorem{lem}[thm]{Lemma}
\newtheorem{cor}[thm]{Corollary}

\theoremstyle{definition}
\newtheorem{defn}[thm]{Definition}

\newtheorem*{conv}{Convention}

\theoremstyle{remark}
\newtheorem{remark}[thm]{Remark}

\numberwithin{equation}{section}

\begin{document}

\title{Sharp metric obstructions for quasi-Einstein metrics}
\author{Jeffrey S. Case}
\thanks{Partially supported by NSF-DMS Grant No.\ 1004394}
\address{Department of Mathematics \\ Princeton University \\ Princeton, NJ 08544}
\email{jscase@math.princeton.edu}
% \date{}
\begin{abstract}
Using the tractor calculus to study smooth metric measure spaces, we adapt results of Gover and Nurowski to give sharp metric obstructions to the existence of quasi-Einstein metrics on suitably generic manifolds.  We do this by introducing an analogue of the Weyl tractor $W$ to the setting of smooth metric measure spaces.  The obstructions we obtain can be realized as tensorial invariants which are polynomial in the Riemann curvature tensor and its divergence.  By taking suitable limits of their tensorial forms, we then find obstructions to the existence of static potentials, generalizing to higher dimensions a result of Bartnik and Tod, and to the existence of potentials for gradient Ricci solitons.
\end{abstract}
\maketitle

%%%%%%%%%%%%%%%%%%%%%%%%%%%%%%%%%%%%%%%%%%%%%%%%%%%%%%%%%%%%%
%                                                           %
% Structure of the document                                 %
%                                                           %
% 1. Intro                                                  %
% *. Acknowledgments                                        %
% 2. Background                                             %
%   a. Tractor Bundles                                      %
%   b. Quasi-Einstein Metrics                               %
% 3. A New Curvature Tractor                                %
% 4. Sharp Tractor Obstructions                             %
% 5. Sharp Tensorial Obstructions                           %
%   a. Weakly generic metrics in pvs Thm                    %
%   b. Static metrics                                       %
%   c. Gradient Ricci solitons                              %
% 5. Concluding remarks                                     %
%                                                           %
%%%%%%%%%%%%%%%%%%%%%%%%%%%%%%%%%%%%%%%%%%%%%%%%%%%%%%%%%%%%%

\section{Introduction}
\label{sec:intro}

Recently, the study of quasi-Einstein metrics has emerged as a unified framework for studying a variety of important metrics in Riemannian geometry.  Roughly speaking, these are metrics $g$ on a manifold such that there exists a function $\phi$ and constants $m\in\bR\cup\{\pm\infty\}$, $\lambda\in\bR$ such that
\begin{equation}
\label{eqn:qe_rough}
\Ric + \nabla^2\phi - \frac{1}{m} d\phi\otimes d\phi = \lambda g ,
\end{equation}
where $\Ric$ is the Ricci curvature of $g$ and $\nabla^2\phi$ is the Hessian of $\phi$.  Depending on the choice of $m$, these include conformally Einstein metrics, gradient Ricci solitons, static metrics in general relativity, and, more generally, the bases of product manifolds which carry a conformally Einstein warped product metric.  Since~\eqref{eqn:qe_rough} gives a family of second-order (degenerate) elliptic PDEs in $g$ and $\phi$ which vary only in the lower order terms, it is natural to try to extend our understanding of certain well-studied quasi-Einstein metrics to all quasi-Einstein metrics.

One natural question to ask is whether or not there are \emph{sharp} obstructions to a given manifold admitting a quasi-Einstein metric, in the sense that the obstruction vanishes if and only if such a metric exists.  As a general question, this is very difficult, with satisfactory answers only known for dimensions $n\leq 2$ in general (see~\cite{Besse,HePetersenWylie2010}) or $n\leq 3$ for Einstein metrics in a given conformal class (see~\cite{Besse}).  Noting that~\eqref{eqn:qe_rough} depends both on a metric and a function, however, one might first simplify the question and ask if there are \emph{metric} obstructions to a given manifold admitting a quasi-Einstein metric; that is, given a fixed metric $g$, can we specify necessary and sufficient conditions for there to exist a function $\phi$ such that~\eqref{eqn:qe_rough} holds?

Within a physical context, this question has attracted a great deal of attention, and there are a number of rather general results for the special cases $m=2-n$ and $m=1$, corresponding to finding an Einstein metric in a conformal class or a static potential in general relativity, respectively.  In the former case, it is well-known that the Cotton-York tensor solves the problem in dimension three, while there are a number of results for suitably generic metrics in dimensions $n\geq 4$ (cf.\ \cite{GoverNurowski2006,KNT1985,Listing2006,Szekeres1963}).  On the other hand, again for suitably generic metrics, Bartnik and Tod~\cite{BartnikTod2006} have solved the problem in dimension three, with the solution similar in many ways to~\cite{KNT1985}.

In the case of conformally Einstein metrics, the work of Gover and Nurowski~\cite{GoverNurowski2006} is the most general.  Their work is in many ways a natural outgrowth of the work of Szekeres~\cite{Szekeres1963}, who studied the four-dimensional case using twistor methods.  In particular, Gover and Nurowsky use the tractor calculus to find obstructions in all dimensions, generalizing those found using tensorial methods by Kozameh, Newman and Tod~\cite{KNT1985} and by Listing~\cite{Listing2006}.

The tractor calculus is an important computational and conceptual tool in conformal geometry, where it generalizes Penrose's twistor theory and gives an intrinsic description of many aspects of the ambient metric of Fefferman and Graham, and thus is a natural analogue of the tensor calculus as used in Riemannian geometry (cf.\ \cite{CapGover2003,FeffermanGraham1985,Penrose1984}).  In particular, conformally Einstein metrics correspond to parallel sections of the so-called standard tractor bundle, whence derivatives of the curvature of this connection yield a number of obstructions to Einstein metrics within a conformal class.  The important insight of Gover and Nurowski~\cite{GoverNurowski2006} is that, subject to a genericity assumption on the metric, one only needs to consider at most one derivative of the curvature to yield sharp obstructions.

The weakest assumption needed in~\cite{GoverNurowski2006} is that the Weyl curvature $W\colon TM\to\Lambda^2T^\ast M\otimes TM$ is injective, which is an open condition on the space of conformal classes on $M$ (hence ``generic'').  From the tractor point of view, the Weyl curvature is the projecting part of the normal tractor connection.  Roughly speaking, this means that for conformal classes which are not generic, there are tractor-valued objects preserved by the conformal holonomy group.  This situation has been considered by Alt~\cite{Alt2006}, where Einstein metrics in certain non-generic conformal classes are found using holonomy methods.

Returning now to quasi-Einstein metrics, the author has recently demonstrated a way to study quasi-Einstein metrics in a unified way using ideas from conformal geometry and by taking suitable ``limits'' as $\lv m\rv\to\infty$~\cite{Case2010a,Case2010b}.  Moreover, when $m$ is finite and not one of three exceptional values, one can define a natural connection and codimension one subbundle of the standard tractor bundle such that parallel sections of the subbundle correspond to conformally quasi-Einstein metrics~\cite{Case2011t}.  Thus, it becomes natural to try to adapt the methods of Gover and Nurowski~\cite{GoverNurowski2006} to find similar sharp metric obstructions to the existence of quasi-Einstein metrics.

The purpose of this article is to carry out this program, and also to consider what happens in the limiting cases $\lv m\rv=\infty$ and $m=1-n$, which do not presently fit into the tractor picture mentioned above.  Analogous to~\cite{GoverNurowski2006}, we do this by introducing the ``weighted Weyl tractor'' associated to a smooth metric measure space (cf.\ \cite{Gover2001}).  Besides being the mechanism by which we find our obstructions, this tensor also seems to suggest another interesting link between conformal geometry and aspects of the Ricci flow via Hamilton's matrix Harnack inequality~\cite{Hamilton1993}, as we will discuss in Section~\ref{sec:conclusion}.

The obstructions constructed using the weighted Weyl tractor are objects which live in tractor bundles, and as such are associated to a conformal class.  However, by fixing a metric in that conformal class, these obstructions can be written in purely tensorial terms, yielding generalizations of the obstructions found in~\cite{GoverNurowski2006,KNT1985,Listing2006}.  Moreover, these tensorial obstructions will admit analogues in the limiting cases $\lv m\rv=\infty$ and $m=1-n$ which cannot be treated using tractor methods, and in particular the latter will yield higher dimensional generalizations of the obstructions to the existence of static metrics found by Bartnik and Tod~\cite{BartnikTod2006} in dimension three.

Finally, we note that while we do not consider the corresponding holonomy groups here, they do exist (see~\cite{Case2011t}).  In particular, it is likely that, similar to~\cite{Alt2006}, one can understand the existence question for quasi-Einstein metrics in many nongeneric conformal classes using holonomy methods.

This article is organized as follows.  In Section~\ref{sec:bg}, we review some essential definitions and ideas from the tractor calculus and the study of quasi-Einstein metrics.  In Section~\ref{sec:curv}, we define the weighted Weyl tractor and use it in Section~\ref{sec:result} to give sharp obstructions to the existence of quasi-Einstein metrics in the tractor language.  In Section~\ref{sec:appl}, we translate the results of Section~\ref{sec:result} into purely tensorial language, and use it to discuss in particular sharp metric obstructions to the existence of static metrics and gradient Ricci solitons.  In Section~\ref{sec:conclusion}, we discuss the relationship between the positivity of the weighted Weyl tractor and Hamilton's matrix Harnack inequality~\cite{Hamilton1993}.

% Acknowledgments
\subsection*{Acknowledgments}
I would like to thank Robert Bartnik and Rod Gover for introducing me to this question, and explaining their respective works.  I would also like to thank the referee for their many useful comments.

\section{Background}
\label{sec:bg}

\subsection{Tractor bundles}
\label{sec:bg/tractor}

We begin by recalling some standard definitions for the tractor calculus in conformal geometry.  We will only take the point of view of vector bundles, following~\cite{Bailey1994,BaumJuhl2010,Case2011t}.  For us, the advantage of using tractor bundles is that they naturally lead to the desired obstructions.  In particular, we will consider the standard tractor bundle $\mT$, the adjoint tractor bundle $\mA=\mT\wedge\mT$, and higher exterior powers of $\mT$, the so-called $k$-form tractors.

A conformal manifold is a pair $(M^n,c)$ of a smooth manifold $M$ together with an equivalence class $c$ of metrics using the equivalence relation $g\sim h$ if and only if $g=e^{2\sigma}h$ for some $\sigma\in C^\infty(M)$.  The conformal class $c$ can be regarded as a ray subbundle $\mC\subset S^2T^\ast M$, sections of which are metrics $g\in c$.  For any $s\in\bR^+$, there is a natural action $\delta_s\colon\mC\to\mC$, defined by $\delta_s(x,g_x)=(x,s^2g_x)$, making $\mC$ into a $\bR^+$-principle bundle.  Given $w\in\bR$, the \emph{conformal density bundle of weight $w$} is the line bundle $E[w]$ associated to $\mC$ via the representation $\bR^+\ni t\mapsto t^{-w/2}\in\End(\bR)$ of $\bR^+$ by $\bR$.

As we will use them, the conformal density bundles $E[w]$ are trivial line bundles (as vector bundles over $M$) associated to the conformal class $c$ with the property that a \emph{choice of scale} $g\in c$ induces an isomorphism $\Gamma\left(E[w]\right)\cong_g C^\infty(M)$ with the property that if $\sigma\in\Gamma\left(E[w]\right)$ and $\sigma_g\in C^\infty(M)$ is the corresponding smooth function under this isomorphism, and if $s\in C^\infty(M)$ determines another metric $h=e^{2s}g\in C^\infty(M)$, then
\begin{equation}
\label{eqn:density_transformation}
\sigma_h=e^{ws}\sigma_g .
\end{equation}
In the interests of brevity, we will henceforth denote by $\mE[w]=\Gamma\left(E[w]\right)$ the space of smooth sections of $E[w]$, and will use the symbol $V$ to denote both a vector bundle over $M$ and its space of smooth functions; in particular, we use the symbol $\mE[w]$ to denote both $E[w]$ and its space of smooth sections.  Additionally, we will more briefly denote the transformation rule~\eqref{eqn:density_transformation} by $\sigma\mapsto e^{ws}\sigma$, with the symbol ``$\mapsto$'' specifying how the function defined by $\sigma$ transforms when one changes metrics $g\mapsto e^{2s}g$.

Given a vector bundle $V$ over $M$, we will frequently denote the corresponding \emph{vector density bundle} $V[w]=V\otimes\mE[w]$.  In this way, one easily sees that the conformal class $c$ determines the \emph{conformal metric} $\bg\in S^2T^\ast M[2]$ by $\bg(\cdot,g):=g$.  This section determines the canonical isomorphism $TM[w]\cong T^\ast M[w+2]$, which we shall use without comment.

Where densities in conformal geometry play a role similar to that of functions in Riemannian geometry, (standard) tractors play a role similar to that of vector fields.  Tractors are sections of the \emph{standard tractor bundle} $\mT$, which can be described as the vector bundle for which a choice of scale $g\in c$ induces an isomorphism $\mT\cong_g\bR\oplus TM\oplus\bR$, together with the transformation law
\begin{equation}
\label{eqn:tractor_transformation}
\begin{pmatrix} \sigma\\\omega\\\rho \end{pmatrix} \mapsto \begin{pmatrix} e^s\sigma\\e^{-s}(\omega+\sigma\,\nabla s)\\e^{-s}(\rho-g(\nabla s,\omega) - \frac{1}{2}|\nabla s|^2\,\sigma) \end{pmatrix} ,
\end{equation}
where the norms and gradients on the right hand side are taken with respect to $g$.  It is clear that given a tractor, the topmost nonvanishing component is conformally invariant, but that the lower components are not; in other words, $\mT$ is a decomposable vector bundle over $M$ associated to $\mC$, but there is no canonical splitting of $\mT$ in terms of the conformal class $c$.  A convenient way to express this in terms of the conformal class $c$ is to regard the standard tractor bundle as a filtered vector bundle
\begin{equation}
\label{eqn:tractor_filtered}
\mT = \mT^{-1} \supset \mT^0 \supset \mT^1,
\end{equation}
where the quotients $\mT^{-1}/\mT^0$, $\mT^0/\mT^1$, and $\mT^1$ correspond to the top, middle, and bottom components, respectively.  On the other hand, \eqref{eqn:tractor_transformation} allows us to easily identify these quotients as density-valued tensor bundles; specifically, $\mT^{-1}/\mT^0=\mE[1]$, $\mT^0/\mT^1=TM[-1]$, and $\mT^1=\mE[-1]$.  We shall write this identification and the filtration~\eqref{eqn:tractor_filtered} as the composition series
\begin{equation}
\label{eqn:tractor_series}
\mT = \mE[1] \semiplus TM[-1] \semiplus \mE[-1] .
\end{equation}
In particular, there is a canonical inclusion $\mE[-1]\hookrightarrow\mT$ and a canonical projection $\mT\to\mE[1]$.

\begin{conv}
Unless otherwise specified, whenever a tractor $I\in\mT$ and a choice of scale $g\in c$ are given, we shall denote the components of $I$ by
\[ I = \begin{pmatrix}\sigma\\\omega\\\rho\end{pmatrix} . \]
\end{conv}

The standard tractor bundle carries a canonical metric $h\in\mT^\ast\otimes\mT^\ast$ and a family of \emph{preferred connections} $\nabla^\prime\colon\mT\to T^\ast M\otimes\mT$, given in a choice of scale by
\begin{align}
\label{eqn:tractor_metric} h(I,I) & = 2\sigma\rho + |\omega|_g^2 \\
\label{eqn:tractor_connection} \nabla^\prime I & = \begin{pmatrix}\nabla\sigma-g(\omega,\cdot)\\\nabla\omega+\sigma P^\prime+\rho\,g\\\nabla\rho-P^\prime(\omega)\end{pmatrix} ,
\end{align}
where $P^\prime$ is any symmetric $(0,2)$-tensor which transforms like the Schouten tensor
\[ P=\frac{1}{n-2}\left(\Ric-\frac{R}{2(n-1)}g\right) ; \]
that is, $P_{\hat g}^\prime-P_{\hat g}=P_g^\prime-P_g$ for all $g,\hat g\in c$.  With such a choice of $P^\prime$, it is easily checked that $\nabla^\prime$ is well-defined; i.e.\ the formula~\eqref{eqn:tractor_connection} transforms according to~\eqref{eqn:tractor_transformation}.  Likewise, one easily checks that the tractor metric is well-defined.

It is easily checked that $\nabla^\prime h=0$ for any preferred connection, and thus may freely identify $\mT\cong\mT^\ast$ using the tractor metric, as we shall do without further comment.  When we are considering the case where $P^\prime=P$, we call the corresponding preferred connection $\nabla$ the \emph{normal tractor connection}, as it corresponds to the normal Cartan connection (cf.\ \cite{CapSlovak2009}).

\begin{conv}
Given a choice of scale $g\in c$, we will sometimes denote a tractor $I\in\mT$ by the triple $I=(\rho,\omega,\sigma)$, according to the isomorphism $\mT\cong\mT^\ast$.  Note, however, that decomposition should be read in the ``opposite order'' of the isomorphism $\mT\cong_g\bR\oplus TM\oplus\bR$; i.e.\ one should regard $(\rho,\omega,\sigma)\in\mE[-1]\oplus TM[-1]\oplus\mE[1]$.

Additionally, we shall usually denote the tractor metric by $\lp\cdot,\cdot\rp$, and the corresponding norm by $\lv\cdot\rv$, and in this way allow ourselves the possibility to use $h$ to describe metrics on other bundles and manifolds.
\end{conv}

The inclusion $\mE[-1]\hookrightarrow\mT$ and the projection $\mT\to\mE[1]$ are both given by a canonical tractor $X\in\mT[1]$, which we shall call the \emph{projector}.  More precisely, the maps
\begin{align*}
I & \mapsto \lp I,X\rp \\
\rho & \mapsto X\otimes\rho
\end{align*}
give the inclusion and the projection, respectively.  In particular, we see that we can canonically write $X=(1,0,0)$ given any choice of scale.

For certain computations and constructions, it will also be necessary to introduce two scale-dependent ``tractors'' which recover the decomposition $\mT\cong_g\mE[1]\oplus TM[-1]\oplus\mE[-1]$.  These are $Y_g\in\mT[-1]$, $Z_g\in T^\ast M\otimes\mT[1]$, and are given by $Y_g=(0,0,1)$, $Z_g(\omega)=(0,\omega,0)$.  When either the scale is clear from context, or $Y$ and $Z$ appear in a scale-independent manner, we shall omit the subscript $g$.

The \emph{tractor-$D$ operator} $\bD\colon\mE[w]\to\mT[w-1]$ defined by
\[ \bD\sigma = \begin{pmatrix} w(n+2w-2)\sigma\\(n+2w-2)\nabla\sigma\\-(\Delta\sigma+\frac{R}{2(n-1)}\sigma) \end{pmatrix} \]
gives a natural way to construct tractors from densities --- this operator is natural in the sense that it is a differential splitting operator for a BGG sequence in conformal geometry; see~\cite{Hammerl2008} and references therein.  Moreover, by replacing the gradient and the Laplacian with the normal tractor connection and the corresponding Laplacian $\Delta=\nabla^\ast\nabla$, one yields the tractor-$D$ operator $\bD\colon\mV[w]\to\mT\otimes\mV[w-1]$ by the same formula, where $\mV$ is any tractor bundle; that is, any combination of tensor products, symmetric tensor products, and skew-symmetric tensor products of $\mT$.  Note that this latter point of view doesn't give rise to a canonical operator.  Similarly, while one can easily define tractor-$D$ operators for preferred connections, the ``best'' definition is not clear without further context.  We will return to this question for our setting in Section~\ref{sec:bg/qe} and Section~\ref{sec:result}.

The algebraic structure underlying the tractor calculus is essentially contained within the \emph{adjoint bundle} $\mA=\mT\wedge\mT$ (cf.\ \cite{CapGover2002}).  The composition series~\eqref{eqn:tractor_series} for the standard tractor bundle yields the composition series
\begin{equation}
\label{eqn:adjoint_series}
\mT\wedge\mT=TM[0]\semiplus\left(\Lambda^2T^\ast M[2]\oplus\mE[0]\right)\semiplus T^\ast M[0] .
\end{equation}
$\mA$ inherits a metric and a family of preferred connections from those on $\mT$.  $\mA$ also contains important algebraic information, via the action $\mA\otimes\mT\mapsto\mT$ given by contraction, $A\otimes I\mapsto A(I):=I\contr A$, as well as an algebraic bracket
\[ \{A,B\}(I,J) = \lp A(I), B(J)\rp - \lp A(J), B(I) \rp . \]
In this way, one readily verifies that fibers $\mA_x$ of the adjoint tractor bundle are isomorphic as Lie algebras to $\kso(n+1,1)$.

The adjoint tractor bundle appears naturally when considering the curvature of a preferred connection $\nabla^\prime$ on $\mT$.  More precisely, its curvature $\mR^\prime$ will be a section of $\Lambda^2T^\ast M\otimes\mA$, and is given in the obvious ``matrix notation'' induced by the definition $\mA=\mT\wedge\mT$ by
\begin{equation}
\label{eqn:curvature_matrix}
\mR^\prime(x,y) = \begin{pmatrix} 0&0&0\\-dP^\prime(x,y)&A^\prime(x,y)&0\\0&dP^\prime(x,y)&0\end{pmatrix} ,
\end{equation}
where $A^\prime=\Rm-P^\prime\wedge g$ (cf.\ \eqref{eqn:curvature_vector}).

\begin{remark}
\label{rk:conventions}
Note that in~\eqref{eqn:curvature_matrix}, we are defining all curvatures using the convention
\[ \mR(x,y) = -\nabla_x\nabla_y+\nabla_y\nabla_x+\nabla_{[x,y]} . \]
With this convention, $\Rm(x,y)w=\imath_w\Rm(x,y)$, and in particular, the curvature $\mR^\prime(I)$ of a tractor is given by the action of $\mA$ on $\mT$ described above.  Note that our convention is the opposite of that found in many references which utilize the tractor calculus (e.g.\ \cite{Bailey1994,GoverNurowski2006}). 
\end{remark}

We can also consider higher exterior powers of $\mT$, yielding \emph{$k$-form tractors} as sections of $\Lambda^k\mT$; in particular, adjoint tractors are simply $2$-form tractors.  Again~\eqref{eqn:tractor_series} induces a composition series
\begin{equation}
\label{eqn:k-form_series}
\Lambda^k\mT = \Lambda^{k-1}T^\ast M[k] \semiplus \left(\Lambda^kT^\ast M[k]\oplus\Lambda^{k-2}T^\ast M[k-2]\right) \semiplus \Lambda^{k-1}T^\ast M[k-2] .
\end{equation}
This in particular will allow us to use a choice of scale $g\in c$ to view a $k$-form tractor $A$ in ``vector notation''
\begin{equation}
\label{eqn:k-form_vector}
A = \begin{pmatrix} \alpha\\\Phi\quad|\quad\phi\\\beta \end{pmatrix} \in \begin{pmatrix} \Lambda^{k-1}T^\ast M[k]\\\Lambda^kT^\ast M[k]\quad|\quad\Lambda^{k-2}T^\ast M[k-2]\\\Lambda^{k-1}T^\ast M[k-2]\end{pmatrix} .
\end{equation}

\begin{remark}
Our convention in defining $A=(\beta,(\Phi,\phi),\alpha)$ will be such that~\eqref{eqn:curvature_matrix} is written in the form~\eqref{eqn:k-form_vector} as
\begin{equation}
\label{eqn:curvature_vector}
\mR^\prime = \begin{pmatrix}0\\A^\prime\quad|\quad0\\-dP^\prime\end{pmatrix} .
\end{equation}
This convention is completely determined by~\eqref{eqn:k-form_action} below.
\end{remark}

The composition series~\eqref{eqn:k-form_series} yields an inclusion $\Lambda^{k-1}T^\ast M[k-2]\hookrightarrow\Lambda^k\mT$ and a projection $\Lambda^k\mT\to\Lambda^{k-1}T^\ast M[k]$, both of which are determined by the canonical section $X\wedge(\wedge^{k-1}Z)\in\Lambda^{k-1}T^\ast M \otimes\Lambda^k\mT[k]$.

It will be useful to have a vector expression for the connection $\nabla^\prime\colon\Lambda^k\mT\to T^\ast M\otimes\Lambda^k\mT$ induced by a preferred connection on $\mT$.  Our conventions are such that
\begin{equation}
\label{eqn:k-form_connection}
\nabla^\prime\begin{pmatrix} \alpha\\\Phi\quad|\quad\phi\\\beta\end{pmatrix} = \begin{pmatrix} \nabla\alpha-(\Phi-g\wedge\phi)\\\nabla\Phi + P^\prime\wedge\alpha + g\wedge\beta\quad|\quad\nabla\phi-\imath_{P^\prime(\cdot)}\alpha + \beta\\\nabla\beta - \imath_{P^\prime(\cdot)}(\Phi-g\wedge\phi)\end{pmatrix} .
\end{equation}

There is also a natural contraction map $\mT\otimes\Lambda^k\mT\to\Lambda^{k-1}\mT$, given by
\begin{subequations}
\label{eqn:k-form_action}
\begin{align}
\imath_X\begin{pmatrix} \alpha\\\Phi\quad|\quad\phi\\\beta\end{pmatrix} & = \begin{pmatrix} 0\\\alpha\quad|\quad0\\\phi\end{pmatrix} \\
\imath_{Z(\omega)}\begin{pmatrix} \alpha\\\Phi\quad|\quad\phi\\\beta\end{pmatrix} & = \begin{pmatrix} -\imath_\omega\alpha\\\imath_\omega\Phi\quad|\quad\imath_\omega\phi\\-\imath_\omega\beta\end{pmatrix} \\
\imath_Y\begin{pmatrix} \alpha\\\Phi\quad|\quad\phi\\\beta\end{pmatrix} & = \begin{pmatrix} -\phi\\\beta\quad|\quad0\\0\end{pmatrix} .
\end{align}
\end{subequations}
Writing out in detail the action of $A=I_1\wedge\dotso\wedge I_k\in\Lambda^k\mT$ on tractors and relating it to~\eqref{eqn:k-form_action} will give explicit formulae for the vector components of $A$, if one is so inclined.  However, we will not need this here.

Finally, the tractor metric induces a metric on the space of $k$-form tractors.  In the vector notation~\eqref{eqn:k-form_vector}, the corresponding norm is given by
\begin{equation}
\label{eqn:k-form_metric}
\left| \begin{pmatrix} \alpha\\\Phi\quad|\quad\phi\\\beta\end{pmatrix} \right|^2 = \lv\Phi\rv^2 - \lv\phi\rv^2 + 2\lp\alpha,\beta\rp .
\end{equation}
Of course, the metric can be explicitly written down by polarizing~\eqref{eqn:k-form_metric}.
\subsection{Quasi-Einstein metrics}
\label{sec:bg/qe}

Let us now briefly review how one regards quasi-Einstein metrics~\eqref{eqn:qe_rough} as distinguished smooth metric measure spaces, and then recall the essential aspects of the tractor formulation of quasi-Einstein metrics given in~\cite{Case2011t}.

\begin{defn}
A \emph{smooth metric measure space (SMMS)} is a four-tuple $(M^n,g,v^m\dvol_g,m)$ of a Riemannian manifold $(M^n,g)$ with its Riemannian volume element $\dvol_g$, a positive smooth function $v\in C^\infty(M)$, and a dimensional parameter $m\in\bR\cup\{\pm\infty\}$.
\end{defn}

There are two important comments to make about this definition.  First, in the case $m=0$, the measure $v^0\dvol_g$ is by convention the usual Riemannian volume element $\dvol_g$, so that we are not adding any additional data to $(M^n,g)$.  Second, when $\lv m\rv<\infty$, we can define $\phi\in C^\infty(M)$ by $v^m=e^{-\phi}$, and it is in this way that one makes sense of a SMMS with $\lv m\rv=\infty$: It is a four-tuple $(M^n,g,e^{-\phi}\dvol_g,\pm\infty)$.  Throughout this article, when the constant $m$ is clear by context, we will simply denote a SMMS as a triple $(M^n,g,v^m\dvol_g)$.  Furthermore, when a SMMS $(M^n,g,v^m\dvol_g)$ has been given, we will always understand the function $\phi$ to be defined by $v^m=e^{-\phi}$, and will frequently use the symbol $\phi$ to denote ``weighted'' objects; i.e.\ geometric objects defined in terms of a SMMS.

The role of the dimensional parameter $m$ is to specify that the measure $v^m\dvol_g$ should be thought of as an $(m+n)$-dimensional measure.  This philosophy guides the definitions of curvature and conformal transformations for SMMS, and can be motivated by considering collapsing sequences of warped product metrics: Given a SMMS $(M^n,g,v^m\dvol_g)$ with $m\in\bN$, choose some Riemannian manifold $(F^m,h)$.  For each $\varepsilon>0$, consider the warped product
\[ \left( M^n\times F^m, \og=g\oplus (\varepsilon v)^2 h \right) . \]
As $\varepsilon\to0$, it is clear that for any $(p,q)\in M\times F$, the pointed metric space $(M^n\times F^m,\og,(p,q))$ converges in the Gromov-Hausdorff sense as $\varepsilon\to0$ to the pointed metric space $(M^n,g,p)$.  However, more is true: In the sense of Cheeger-Colding~\cite{CheegerColding1997}, we have that
\[ \left( M^n\times F^m, \og, \overline{\dvol}, (p,q)\right) \to \left( M^n,g,v^m\dvol_g,p \right) \]
as $\varepsilon\to0$, where $\overline{\dvol}$ denotes the renormalized Riemannian volume element of $(M\times F,\og)$ (see~\cite{CheegerColding1997} for details).  In this way, we can think of the notion of a SMMS as a generalization of these types of collapsed limits, where we allow the total dimension of $M\times F$ to be any real number, positive or negative, or possibly even infinite.  Moreover, when $m\in\bN$, most of the geometric objects we define can be realized from the horizontal components (i.e.\ the parts depending only on $M$) of their corresponding geometric objects in the warped product $(M\times F^m,\og)$.

The most natural definition associated to a SMMS is that of the weighted divergence.

\begin{defn}
\label{defn:weighted_divergence}
Let $(M^n,g,v^m\dvol_g)$ be a SMMS.  Given a vector bundle $(E,h,\nabla)$ over $M$, let $d\colon\Lambda^kT^\ast M\otimes E\to\Lambda^{k+1}T^\ast M\otimes E$ be the usual twisted exterior derivative.  The \emph{weighted divergence $\delta_\phi\colon\Lambda^{k+1}T^\ast M\otimes E\to\Lambda^kT^\ast M\otimes E$} is defined by
\[ \int_M \lp \xi, d\zeta\rp e^{-\phi}\dvol = -\int_M \lp \delta_\phi\xi,\zeta\rp e^{-\phi}\dvol \]
for all $\xi\in\Lambda^{k+1}T^\ast M\otimes E$ and $\zeta\in\Lambda^kT^\ast M\otimes E$, at least one of which is compactly supported, where $\lp\cdot,\cdot\rp$ denotes the natural pointwise inner product induced by $g$ and $h$ on the corresponding vector bundle.

The \emph{weighted Laplacian $\Delta_\phi$} is the operator $\Delta_\phi\colon C^\infty(M)\to C^\infty(M)$ defined by $\Delta_\phi=\delta_\phi d$.
\end{defn}

In particular, given $u\in C^\infty(M)$, we can write $\Delta_\phi u=\Delta u - \lp\nabla\phi,\nabla u\rp$ in terms of the usual Laplacian of the Riemannian manifold $(M^n,g)$.  Note also that the definition of the weighted divergence makes sense for more general differential operators; however, in this article we shall only need its definition as the (negative of the) adjoint of the twisted exterior derivative.

To study the geometry of a SMMS, one needs useful notions of curvature.  The most common such notion on a SMMS is the Bakry-\'Emery Ricci tensor, which is the weighted analogue of the Ricci curvature and is typically derived as the curvature term appearing in the Bochner (in)equality involving $\Delta_\phi\lv\nabla u\rv^2$ (cf.\ \cite{BakryEmery1985}).

\begin{defn}
Let $(M^n,g,v^m\dvol_g)$ be a SMMS.  The \emph{Bakry-\'Emery Ricci tensor $\Ric_\phi^m$} and the \emph{weighted scalar curvature $R_\phi^m$} are defined by
\begin{align*}
\Ric_\phi^m & = \Ric - mv^{-1}\nabla^2v = \Ric + \nabla^2\phi - \frac{1}{m}d\phi\otimes d\phi \\
R_\phi^m & = R-2mv^{-1}\Delta v-m(m-1)v^{-2}\lv\nabla v\rv^2 = R + 2\Delta\phi - \frac{m+1}{m}\lv\nabla\phi\rv^2 .
\end{align*}

We say that $(M^n,g,v^m\dvol_g)$ is \emph{quasi-Einstein} if there is a constant $\lambda\in\bR$ such that
\begin{equation}
\label{eqn:qe_defn}
\Ric_\phi^m = \lambda g ,
\end{equation}
and call $\lambda$ the \emph{quasi-Einstein constant}.

If $(M^n,g,v^m\dvol_g)$ is quasi-Einstein with quasi-Einstein constant $\lambda$, the \emph{\charconstant\ $\mu$} is the constant such that
\begin{equation}
\label{eqn:charconst_defn}
R_\phi^m + m\mu v^{-2} = (m+n)\lambda .
\end{equation}
\end{defn}

\begin{remark}
When $\lv m\rv=\infty$, by convention the \charconstant\ is the quasi-Einstein constant.
\end{remark}

Given a quasi-Einstein SMMS $(M^n,g,v^m\dvol)$ with quasi-Einstein constant $\lambda$, the existence of a constant $\mu$ such that~\eqref{eqn:charconst_defn} holds is a result of D.-S.\ Kim and Y.\ H.\ Kim~\cite{Kim_Kim}.

It is important to note that, in general, $R_\phi^m\not=\tr_g\Ric_\phi^m$.  Instead, if $m\in\bN$ and one computes in the warped product $(M\times_v F^m,\og=g\oplus v^2h)$, one finds that the horizontal component of $\Ric_{\og}$ is exactly $\Ric_\phi^m$, while if one assumes additionally that $h$ is scalar flat, then the scalar curvature $R_{\og}$ of $\og$ is exactly $R_\phi^m$.  In other words, $R_\phi^m$ includes in addition to the trace of $\Ric_\phi^m$ some information which is ``hidden'' in the fibers $(F,h)$.  Moreover, if $(M^n,g,v^m\dvol)$ is quasi-Einstein with \charconstant\ $\mu$ and $m\in\bN$, whenever $(F^m,h)$ is an Einstein manifold satisfying $\Ric_h=\mu h$, the warped product $(M\times_v F,\og)$ is Einstein with $\Ric_{\og}=\lambda\og$ (cf.\ \cite{Besse,Case_Shu_Wei,Kim_Kim}).

In order to adapt tools from conformal geometry to the study of SMMS, we first must define what it means for two SMMS to be (pointwise) conformally equivalent.

\begin{defn}
Two SMMS $(M^n,g,v^m\dvol_g)$ and $(M^n,\hat g,\hat v^m\dvol_{\hat g})$ are \emph{(pointwise) conformally equivalent} if there is a function $s\in C^\infty(M)$ such that
\begin{equation}
\label{eqn:scms}
\left( M^n,\hat g,\hat v^m\dvol_{\hat g}\right) = \left( M^n, u^{-2}g, u^{-m-n}v^m\dvol_g \right) .
\end{equation}
\end{defn}

There are two observations to make about this definition.  First, \eqref{eqn:scms} is in line with the philosophy that the measure $v^m\dvol_g$ should be regarded as an $(m+n)$-dimensional measure, as it transforms after a conformal change of metric as would the Riemannian volume element of an $(m+n)$-dimensional manifold.  Second, one can regard~\eqref{eqn:scms} as equivalently specifying that a pointwise conformal change on a SMMS is a map $(g,v)\mapsto (e^{2s}g,e^sv)$ for some $s\in C^\infty(M)$.  In this way, we can talk about the conformal class of a SMMS, which yields the definition of a smooth conformal measure space.

\begin{defn}
A \emph{smooth conformal measure space (SCMS)} is a triple $(M^n,c,v^md\nu)$ of a conformal manifold $(M^n,c)$ with $d\nu\in\Lambda^nT^\ast M[n]$ the canonical conformal volume element --- that is, for any choice of scale $g\in c$, $d\nu$ trivializes as the Riemannian volume element $\dvol_g$ --- together with a dimensional parameter $m\in\bR$ and a positive density $v\in\mE[1]$.

A \emph{choice of scale} for $(M^n,c,v^md\nu)$ is a choice of metric $g\in c$.  The induced trivialization of $\mE[1]$ allows us to regard $v\in C^\infty(M)$, in which case the choice of scale $g$ yields the SMMS $(M^n,g,v^m\dvol_g)$.
\end{defn}

One can also regard a choice of scale as a specific choice of function $v\in C^\infty(M)$ represented by the density $v\in\mE[1]$.  We will only ever use this perspective to choose the scale $v=1$, which is equivalent to the choice of scale $v^{-2}\bg\in c$.  Note also that in the definition of a SCMS we are not allowing $\lv m\rv$ to be infinite.  This is because, while one can make sense of~\eqref{eqn:scms} when $\lv m\rv=\infty$, this requires no longer changing the metric conformally (cf.\ \cite{Case2010a}).  Also note that, like the definition of a SMMS, when $m=0$ we are specifying exactly the usual data $(M^n,c,d\nu)$ of a conformal manifold with its canonical volume element; i.e.\ a SCMS with $m=0$ is simply a conformal manifold.

A natural question to ask about a SCMS is whether it admits any quasi-Einstein scales; equivalently, whether there exist $u\in\mE[1]$ such that the choice of scale $u^{-2}\bg\in c$, where defined, gives a quasi-Einstein SMMS.  This question reduces to a conformally invariant PDE in $u$ which is easily computed.

\begin{prop}[\cite{Case2010a}]
Let $(M^n,g,v^m\dvol)$ be a SMMS and suppose that $u\in C^\infty(M)$ is a positive function such that the SMMS defined by~\eqref{eqn:scms} is quasi-Einstein with quasi-Einstein constant $\lambda$ and \charconstant\ $\mu$.  Then
\begin{subequations}
\label{eqn:qe_conf}
\begin{align}
\label{eqn:qe_tf} 0 & = \left(uv\Ric + (m+n-2)v\nabla^2u - mu\nabla^2v\right)_0 \\
\label{eqn:qe_lambda} n\lambda v^2 & = (uv)^2R + (m+2n-2)uv^2\Delta u - mu^2v\Delta v \\
\notag & \quad - (m+n-1)nv^2\lv\nabla u\rv^2 + mnuv\lp\nabla u,\nabla v\rp \\
\label{eqn:qe_mu} n\mu u^2 & = (uv)^2R + (m+n-2)uv^2\Delta u - (m-n)u^2v\Delta v \\
\notag & \quad - (m+n-2)nuv\lp\nabla u,\nabla v\rp + n(m-1)u^2\lv\nabla v\rv^2,
\end{align}
\end{subequations}
where $T_0=T-\frac{1}{n}\tr_g T\,g$ denotes the tracefree part of a symmetric $(0,2)$ tensor and all derivatives are computed with respect to $g$.
\end{prop}

There are a few important observations to make about~\eqref{eqn:qe_conf}.

First, these equations make perfect sense if $u$ is allowed to change sign, and in particular the constants $\lambda$ and $\mu$ are well-defined globally in this case.

Second, it is easily verified that so long as one regards $u\in\mE[1]$, \eqref{eqn:qe_conf} is a conformally invariant system of equations.

Third, the conditions~\eqref{eqn:qe_lambda} and~\eqref{eqn:qe_mu} imply that one can regard $\lambda$ and $\mu$ as the squared lengths of $u$ and $v$, respectively.  More precisely, if one sets $\hat u=cu$ and $\hat v=kv$ for constants $c,k>0$, then~\eqref{eqn:qe_conf} still holds provided one replaces $\lambda$ by $\hat\lambda=c^2\lambda$ and $\mu$ by $\hat\mu=k^2\mu$.  For this reason, when we consider fixing the measure $v^m\dvol_g$, we will frequently also fix the desired \charconstant\ $\mu$.  When we do this, we will say that we have fixed a SMMS $(M^n,g,v^m\dvol)$ with \charconstant\ $\mu$, so that the search for a function $u\in C^\infty(M)$ satisfying~\eqref{eqn:qe_conf} can be rephrased as the search for a quasi-Einstein scale as follows:

\begin{defn}
Let $(M^n,g,v^m\dvol_g)$ be a SMMS with \charconstant\ $\mu$.  Then $u\in C^\infty(M)$ is a \emph{quasi-Einstein scale} if the SMMS
\[ \left( M^n,u^{-2}g,u^{-m-n}v^m\dvol_g\right) \]
is a quasi-Einstein SMMS with \charconstant\ $\mu$ wherever it is defined.

The \emph{quasi-Einstein constant $\lambda$} of such a quasi-Einstein scale is the constant $\lambda$ determined by~\eqref{eqn:qe_conf}.
\end{defn}

Note in particular that this definition allows $u$ to change signs.  Similar terminology will be used for SCMS $(M^n,g,v^md\nu)$.

Finally, there is a symmetry hidden in~\eqref{eqn:qe_conf} which, while still not completely understood, is nevertheless important in applying our tractor methods to the study of quasi-Einstein metrics.

\begin{cor}[\cite{Case2010a}]
\label{cor:duality}
Let $(M^n,g,v^m\dvol_g)$ be a SMMS with \charconstant\ $\mu$ which admits a positive quasi-Einstein scale $u\in C^\infty(M)$ with quasi-Einstein constant $\lambda$.  Then the SMMS $(M^n,g,u^{2-m-n}\dvol_g)$ is a SMMS with \charconstant\ $\lambda$ such that $v$ is a quasi-Einstein scale with quasi-Einstein constant $\mu$.
\end{cor}

The proof is to simply observe that~\eqref{eqn:qe_conf} is invariant under the transformation
\[ (u,v,m,\lambda,\mu) \mapsto (v,u,2-m-n,\mu,\lambda) . \]
The main benefit of the ``duality'' of Corollary~\ref{cor:duality} is that, given a fixed Riemannian manifold $(M^n,g)$ and a fixed dimensional parameter $m\geq 0$, it allows one to study the space of functions $v\in C^\infty(M)$ for which the SMMS $(M^n,g,v^m\dvol_g)$ is quasi-Einstein with quasi-Einstein constant $\lambda$ as the space of quasi-Einstein scales for the SMMS $(M^n,g,1^{2-m-n}\dvol_g)$ with \charconstant\ $\lambda$.  This is precisely the perspective of He, Petersen and Wylie~\cite{HePetersenWylie2010,HePetersenWylie2010b,HePetersenWylie2011c}; see Section~\ref{sec:appl/static} for a useful application of this idea to static metrics.

In order to generalize the work of Gover and Nurowski~\cite{GoverNurowski2006} to the setting of SCMS, we must first describe the relevant tractor objects associated to such a space.  For these definitions to make sense, we must additionally assume that the dimensional parameter $m\not\in\{-n,1-n,2-n\}$, an assumption which we shall make implicitly for the remainder of this article.

\begin{defn}
Let $(M^n,c,v^md\nu)$ be a SCMS with \charconstant\ $\mu$.  The \emph{SCMS-scale tractor $\tilde J$} is defined by
\begin{equation}
\label{eqn:tildeJ}
\tilde J = J + \frac{(m+2n-2)(\mu-(m-1)\lv J\rv^2)}{2(m+n-1)(m+n-2)v}X,
\end{equation}
where $J=\frac{1}{n}\bD v$.

The \emph{$W$-tractor subbundle $\mT^W$} is the codimension one subbundle $\mT^W=\tilde J^\perp\subset\mT$

The \emph{$W$-tractor connection $\nabla^W\colon\mT\to T^\ast M\otimes\mT$} is defined, given a choice of scale $g\in c$, by
\begin{equation}
\label{eqn:nablaW}
\nabla^W\begin{pmatrix}\sigma\\\omega\\\rho\end{pmatrix} = \begin{pmatrix}\nabla\sigma-\omega\\\nabla\omega+\sigma P^W + \rho g\\\nabla\rho - P^W(\omega)\end{pmatrix},
\end{equation}
where the \emph{weighted Schouten tensor $P^W\in S^2T^\ast M$} is defined by
\begin{equation}
\label{eqn:weighted_schouten}
P^W = \frac{1}{m+n-2}\left(\Ric_\phi^m - \J^W g\right), \quad \J^W = \frac{R_\phi^m+m\mu v^{-2}}{2(m+n-1)} .
\end{equation}
We say that a tractor $I\in\mT$ is \emph{$W$-parallel} if $\nabla^WI=0$.

The \emph{$W$-tractor-$D$ operator $\bD^W\colon\mV[w]\to\mT\otimes\mV[w-1]$} is defined, given a choice of scale $g\in c$, by
\begin{equation}
\label{eqn:DW}
\bD^W\sigma = \begin{pmatrix} w(m+n+2w-2)\sigma\\(m+n+2w-2)\nabla^W\sigma\\-(\Delta_\phi^W\sigma+w\J^W\sigma)\end{pmatrix},
\end{equation}
where $\Delta_\phi^W=\delta_\phi^W d^W$ is computed using the twisted exterior derivative $d^W\colon\mV\to T^\ast M\otimes\mV$ defined using $\nabla^W$.
\end{defn}

That $\nabla^W$ and $\bD^W\colon\mE[w]\to\mT[w-1]$ are well-defined is proven in~\cite{Case2011t}.  That $\bD^W$ is well-defined as an operator on tractor bundles as given here is easily checked by direct computation.  Note also that all of the weighted tractor objects defined above depend on the choice of the \charconstant\ $\mu$.

Arguably the main reason for the introduction of these weighted tractor objects is that they give a correspondence between quasi-Einstein scales and parallel tractors.

\begin{thm}
\label{thm:equivalence}
Let $(M^n,c,v^md\nu)$ be a SCMS with \charconstant\ $\mu$.  There is a one-to-one correspondence between quasi-Einstein scales $u\in\mE[1]$ and $W$-parallel tractors $I\in\mT^W$, given by
\[ u \longleftrightarrow \frac{1}{m+n}\bD^W u = I . \]
\end{thm}

In particular, Theorem~\ref{thm:equivalence} states that if $I\in\mT^W$ is $W$-parallel, then there exists a function $u\in\mE[1]$ such that $I=\frac{1}{m+n}\bD^Wu$.  This is closely related to the idea that the $W$-tractor connection is a prolongation of (part of) the conformally quasi-Einstein condition~\eqref{eqn:qe_conf}.  Another useful realization of this idea, which we will find useful in the present article, is the following lemma from~\cite{Case2011t}.

\begin{lem}
\label{lem:cute_algebraic_trick}
Let $(M^n,c,v^md\nu)$ be a SCMS with \charconstant\ $\mu$.  Suppose that $I\in\mT^W$ is such that $\nabla^WI=\beta\otimes X$ for some one-form $\beta$.  Then $\nabla^WI=0$; i.e.\ $\beta=0$.
\end{lem}

Finally, in order to introduce the weighted analogue of the Weyl tractor, the following formulae from~\cite{Case2011t} will be useful.

\begin{lem}
Let $(M^n,g,v^m\dvol)$ be a SMMS with \charconstant\ $\mu$.  Let $P^W$ be the weighted Schouten tensor~\eqref{eqn:weighted_schouten} and define the weighted Weyl curvature $A^W$ by $A^W=\Rm-P^W\wedge g$; i.e.\ $A^W$ is the topmost component of the curvature~\eqref{eqn:curvature_matrix} of the $W$-tractor connection computed in the scale $g$.  Using the scale $g$ to write $\tilde J$ in~\eqref{eqn:tildeJ} as $\tilde J=(\tilde y,\nabla v,v)$, it holds that
\begin{align}
\label{eqn:trdP^C} \tr_g dP^W & = mv^{-1}\left(\nabla\tilde y - P^W(\nabla v)\right) \\
\label{eqn:trA^C} \tr_g A^W & = mv^{-1}\left(vP^W+\nabla^2 v + \tilde y \, g\right) \\
\label{eqn:divA^C} \delta_\phi A^W & = (m+n-3)dP^W - mv^{-2}\nabla v\wedge \left(vP^W+\nabla^2v+\tilde y\,g \right) \\
\label{eqn:divdP^C} \delta_\phi dP^W & = -mv^{-2}\nabla v\wedge \left(\nabla\tilde y - P^W(\nabla v)\right) ,
\end{align}
where we regard~\eqref{eqn:divdP^C} as an expression relating sections of $\Lambda^2T^\ast M$.
\end{lem}
\section{A New Curvature Tractor}
\label{sec:curv}

In order to use the tractor calculus introduced in Section~\ref{sec:bg} to find our sharp obstructions, we will need to introduce a new curvature-like tractor to the study of smooth conformal measure spaces.  As we will see, this tractor is a section of $\mA\otimes\mA[-2]$ which has Riemann tensor symmetries (cf.\ equation~(30) in~\cite{Gover2001}).

\begin{defn}
\label{defn:curvature_tractor}
Let $(M^n,c,v^md\nu)$ be a SCMS with \charconstant\ $\mu$.  The \emph{weighted Weyl tractor $W^W\in\mA\otimes\mA[-2]$} is determined in a choice of scale $g\in c$ by
\begin{equation}
\label{eqn:WC_coords}
W^W = \begin{pmatrix} 0 \\ (m+n-4)\mR^W \quad|\quad 0 \\ -\delta_\phi^W\mR^W + mv^{-2}\nabla^W\tilde J\wedge\tilde J\end{pmatrix},
\end{equation}
where $\delta_\phi^W\colon\Lambda^2T^\ast M\otimes\mA\to T^\ast M\otimes\mA$ is the weighted divergence of the twisted exterior derivative $d^W\colon T^\ast M\otimes\mA\to\Lambda^2T^\ast M\otimes\mA$ defined using the connection $\nabla^W$ on $\mA$, and in~\eqref{eqn:WC_coords} we use the scale $g$ to regard
\[ W^W \in \left( T^\ast M[2] \oplus (\Lambda^2T^\ast M[2]\oplus\mE[0]) \oplus T^\ast M[0]\right) \otimes \mA \]
using~\eqref{eqn:adjoint_series}.
\end{defn}

Observe that in the case $m=0$, corresponding to the familiar setting of conformal geometry, the weighted Weyl tractor is precisely the Weyl tractor $W$ defined in~\cite{Gover2001}.

We claim that $W^W$ is a well-defined section of $\mA\otimes\mA[-2]$ and that $W^W$ is an algebraic curvature tractor, in the sense that it has the usual Riemann curvature symmetries.  The remainder of this section is devoted to establishing these facts.

First, we introduce a useful piece of notation which we shall use in the remainder of this article.  Given sections $A\in\Lambda^2T^\ast M\otimes\Lambda^2T^\ast M$ and $T\in T^\ast M\otimes T^\ast M$ on a Riemannian manifold $(M^n,g)$, there is a natural pairing $\lp A,T\rp\in T^\ast M\otimes T^\ast M$ which arises by contracting in each corresponding factor.  More precisely, given an orthonormal basis $\{e_i\}\subset T_pM$ and vectors $x,y\in T_pM$, we define
\begin{equation}
\label{eqn:natural_pairing}
\lp A,T\rp(x,y) = \sum_{i,j=1}^n A(e_i,x,e_j,y) T(e_i,e_j) .
\end{equation}
Similarly, the tractor metric induces a natural pairing
\[ \left(\mA\otimes\mA\right) \otimes \left(\mT\otimes\mT\right) \to \mT\otimes\mT \]
which we will likewise denote $\lp A,T\rp$ for $A\in\mA\otimes\mA$ and $T\in\mT\otimes\mT$.

Let us first show that $W^W$ is an algebraic curvature tractor in the scale $g\in c$ used in Definition~\ref{defn:curvature_tractor}; that $W^W$ is in fact an algebraic curvature tractor will then follow once we establish that $W^W$ is well-defined.  To see that $W^W$ has Riemann curvature symmetries in the scale $g\in c$, we first find a tensorial formula for the bottom component of~\eqref{eqn:WC_coords}.

\begin{lem}
\label{lem:tractor_formula_ww}
Let $(M^n,c,v^md\nu)$ be a SCMS with \charconstant\ $\mu$.  Given a scale $g\in c$, it holds that
\begin{equation}
\label{eqn:divphi_mR^C}
\delta_\phi^W\mR^W-mv^{-2}\nabla^W\tilde J\wedge\tilde J = \begin{pmatrix} 0\\(m+n-4)dP^W\quad|\quad0\\-B^W\end{pmatrix} ,
\end{equation}
where we have written
\begin{equation}
\label{eqn:B^C}
B^W = \delta_\phi dP^W + v^{-1}\tr dP^W\otimes dv + \lp A^W, P^W-v^{-1}\tilde y g\rp
\end{equation}
for $\delta_\phi dP^W\in T^\ast M\otimes T^\ast M$ via the composition
\[ \Lambda^2 T^\ast M \otimes T^\ast M \xrightarrow{\delta_\phi} T^\ast M\otimes T^\ast M . \]
\end{lem}

\begin{remark}
When $m=0$, so that $W^W$ is the usual Weyl tractor, the tensor $B^W$ defined by~\eqref{eqn:B^C} is the Bach tensor.  This suggests thinking of $B^W$ as the ``weighted Bach tensor,'' though we do not know of any references which have made use of this perspective.
\end{remark}

\begin{proof}

Using the induced connection $\nabla^W$ on $\Lambda^2T^\ast M\otimes\mA$ and~\eqref{eqn:curvature_vector}, we compute that
\begin{equation}
\label{eqn:div_mR^C}
\delta_\phi^W\mR^W = \begin{pmatrix} -\tr A^W\\\quad\delta_\phi A^W-dP^W\quad|\quad-\tr dP^W\\-\delta_\phi dP^W-\lp A^W,P^W\rp \end{pmatrix} .
\end{equation}
On the other hand, \eqref{eqn:trdP^C}, \eqref{eqn:trA^C}, and~\eqref{eqn:divA^C} imply that
\[ -mv^{-2}\nabla^W\tilde J\wedge\tilde J = \begin{pmatrix}\tr A^W\\v^{-1}dv\wedge\tr A^W\quad|\quad\tr dP^W\\v^{-1}\tilde y\tr A^W - v^{-1}\tr dP^W\otimes dv\end{pmatrix}, \]
Combining these facts yields the desired result.
\end{proof}

In order to show that $W^W$ is an algebraic curvature tractor, it remains to show that $B^W$ is symmetric.

\begin{lem}
\label{lem:weighted_bach_symmetric}
Let $(M^n,g,v^m\dvol_g)$ be a SMMS with \charconstant\ $\mu$.  The tensor $B^W$ defined by~\eqref{eqn:B^C} is symmetric.
\end{lem}

\begin{proof}

Since $A^W$ and $P^W$ are symmetric, $\lp A^W,P^W\rp$ is symmetric.  Thus it suffices to show that
\begin{equation}
\label{eqn:bach_goal}
\delta_\phi dP^W + v^{-1}\tr dP^W\otimes dv
\end{equation}
is symmetric.  By the definition of the divergence $\delta\colon\Lambda^2T^\ast M\otimes T^\ast M\to T^\ast M\otimes T^\ast M$,
\begin{equation}
\label{eqn:skew}
\begin{split}
\delta dP^W(x,y) - \delta dP^W(y,x) & = \sum_{i=1}^n \nabla_{e_i}dP^W(e_i,x,y) - \nabla_{e_i} dP^W(e_i,y,x) \\
& = -\sum_{i=1}^n \nabla_{e_i} dP^W(x,y,e_i)
\end{split}
\end{equation}
for any orthonormal basis $\{e_i\}$ of $T_pM$, where the second equality follows from the fact
\begin{equation}
\label{eqn:dP^W_bianchi}
dP^W(x,y,z) + dP^W(y,z,x) + dP^W(z,x,y) = 0
\end{equation}
for all $x,y,z\in TM$.  In particular, observe that~\eqref{eqn:skew} and~\eqref{eqn:dP^W_bianchi} imply that
\[ \delta_\phi dP^W(x,y) - \delta_\phi dP^W(y,x) = -\delta_\phi dP^W \in \Lambda^2T^\ast M, \]
where by an abuse of notation we regard $\delta_\phi dP^W\in T^\ast M\otimes T^\ast M$ on the left hand side and $\delta_\phi dP^W\in\Lambda^2T^\ast M$ on the right hand side as the two different ways one can take the weighted divergence of $dP^W$.  The symmetry of~\eqref{eqn:bach_goal}, and hence of $B^W$, then follows immediately from~\eqref{eqn:divdP^C}.
\end{proof}

\begin{cor}
\label{cor:algebraic_tractor}
Let $(M^n,c,v^md\nu)$ be a SCMS with \charconstant\ $\mu$.  In the scale $g\in c$, the Weyl tractor $W^W$ is an algebraic curvature tractor.
\end{cor}

\begin{proof}

Given $T_1,T_2\in\mA$, define the forms $\alpha_i,\beta_i\in T^\ast M$, $\phi_1\in C^\infty(M)$, and $\Phi_i\in\Lambda^2T^\ast M$ for $i=1,2$ by
\[ T_1 = \begin{pmatrix}\alpha_1\\\Phi_1\quad|\quad\phi_1\\\beta_1\end{pmatrix}, \qquad
   T_2 = \begin{pmatrix}\alpha_2\\\Phi_2\quad|\quad\phi_2\\\beta_2\end{pmatrix} , \]
so that
\begin{equation}
\label{eqn:w_eval}
\begin{split}
W^W\left(T_1,T_2\right) & = (m+n-4)\lp A^W,\Phi_1\otimes\Phi_2\rp + B^W(\alpha_1,\alpha_2) \\
& \quad - (m+n-4)\lp dP^W,\Phi_1\otimes\alpha_2+\Phi_2\otimes\alpha_1\rp ,
\end{split}
\end{equation}
By Lemma~\ref{lem:weighted_bach_symmetric}, $B^W$ is symmetric.  Thus, since $A^W$ is symmetric, it follows that $W^W$ is symmetric.

Using~\eqref{eqn:w_eval} again, it is straightforward to check that given any $I_1,\dotsc,I_4\in\mT$,
\[ W^W(I_1\wedge I_2,I_3\wedge I_4) + W^W(I_2\wedge I_3,I_1\wedge I_4) + W^W(I_3\wedge I_1,I_2\wedge I_4) = 0 . \]
Hence $W^W$ satisfies the Bianchi identity, and is thus an algebraic curvature tractor.
\end{proof}

Let us now turn to the problem of showing that the weighted Weyl tractor is well-defined.  We accomplish this by finding an expression for $W^W$ written purely in terms of tractor operators.

\begin{prop}
\label{prop:well_defined}
Let $(M^n,c,v^md\nu)$ be a SCMS with \charconstant\ $\mu$.  The weighted Weyl tractor can be written
\begin{equation}
\label{eqn:W^W_invariant}
W^W = \frac{1}{m+n-2}\left(\tr\bD^W + m(m+n-3)v^{-1}\tilde J\contr\right)(X\wedge\mR^W) + mv^{-2}X\wedge\nabla^W\tilde J\wedge\tilde J .
\end{equation}
In particular, the weighted Weyl tractor is a well-defined section of $\mA\otimes\mA[-2]$.
\end{prop}

Before we give the proof, let us make some comments about our notation.  From~\eqref{eqn:k-form_series}, we know that there is a well-defined inclusion $\Lambda^{k-1}T^\ast M\hookrightarrow\Lambda^{k}\mT[2-k]$ given by $X\wedge(\wedge^{k-1}Z)\in\Lambda^{k-1}T^\ast M\otimes\Lambda^k\mT[k]$.  This extends to inclusions
\[ \Lambda^2T^\ast M\otimes\mA \hookrightarrow \Lambda^3\mT\otimes\mA[-1], \quad T^\ast M\otimes\mA\hookrightarrow\mA\otimes\mA . \]
In~\eqref{eqn:W^W_invariant}, we are denoting
\begin{align*}
X\wedge\mR^W & \in \Lambda^3\mT\otimes\mA[-1] \\
X\wedge\nabla^W\tilde J\wedge\tilde J & \in \mA\otimes\mA
\end{align*}
for the respective inclusions of $\mR^W\in\Lambda^2T^\ast M\otimes\mA$ and $\nabla^W\tilde J\wedge\tilde J\in T^\ast M\otimes\mA$ in this way.  The contraction
\[ \tilde J\contr (X\wedge\mR^W) \in \mA\otimes\mA[-1] \]
is the usual contraction of a tractor into the first factor, while the operator $\tr\bD^W$ is defined as the composition
\[ \Lambda^3\mT\otimes\mA[-1] \xrightarrow{\bD^W} \mT\otimes\Lambda^3\mT\otimes\mA[-2] \xrightarrow{\tr} \mA\otimes\mA[-2], \]
with the symbol $\tr$ denoting the obvious contraction $\mT\otimes\Lambda^3\mT\mapsto\mA$.

\begin{proof}[Proof of Proposition~\ref{prop:well_defined}]

First observe that, by our conventions described above, it is clear that the right hand side of~\eqref{eqn:W^W_invariant} is a section of $\mA\otimes\mA[-2]$.  To complete the proof, we fix a scale $g\in c$ and compute the right hand side of~\eqref{eqn:W^W_invariant} in the same ``coordinates'' used to express~\eqref{eqn:WC_coords}.

By definition of the inclusion $\Lambda^2T^\ast M\otimes\mA\hookrightarrow\Lambda^3\mT\otimes\mA[-1]$, we have that
\begin{equation}
\label{eqn:xwedger_vector}
\Omega^W := X\wedge\mR^W = \begin{pmatrix}0\\0\quad|\quad0\\\mR^W\end{pmatrix} \in \Lambda^3\mT\otimes\mA[-1] .
\end{equation}
It follows immediately from~\eqref{eqn:k-form_action} that
\begin{equation}
\label{eqn:xwedger_contr}
mv^{-1}\tilde J\contr\Omega^W = \begin{pmatrix} 0\\m\mR^W\quad|\quad0\\-mv^{-1}\imath_{\nabla v}\mR^W\end{pmatrix} \in \mA\otimes\mA[-2] .
\end{equation}

On the other hand, using~\eqref{eqn:k-form_action}, \eqref{eqn:DW}, and the definition of the trace $\tr\colon\mT\otimes\Lambda^3\mT\to\mA$, we see that
\begin{equation}
\label{eqn:xwedger_tr}
\begin{split}
\tr\bD^W\Omega^W & = -(m+n-4)\begin{pmatrix}0\\\mR^W\quad|\quad0\\0\end{pmatrix} - \imath_X \Delta_\phi^W\Omega^W \\
& \quad + (m+n-4)\sum_{i=1}^n \lp \nabla_{e^i}^W\Omega^W, Z(e^i)\rp,
\end{split}
\end{equation}
By the definition of $\nabla^W$ on $\Lambda^3\mT\otimes\mA[-1]$, we see that
\begin{align*}
\nabla_y^W\Omega^W & = \begin{pmatrix} 0\\g(y,\cdot)\wedge\mR^W\quad|\quad\imath_y\mR^W\\\nabla_y^W\mR^W\end{pmatrix} \\
\Delta_\phi^W\Omega^W & = \begin{pmatrix} -(n-4)\mR^W\\\ast\quad|\quad2\delta_\phi^W\mR^W+\imath_{\nabla\phi}\mR^W\\\ast\end{pmatrix} ,
\end{align*}
where the terms marked $\ast$ are irrelevant, as they disappear upon contraction with $X$ in~\eqref{eqn:xwedger_tr}.  In particular, using these to evaluate the right hand side of~\eqref{eqn:xwedger_tr}, it follows that
\begin{equation}
\label{eqn:w_components}
W^W = \begin{pmatrix} 0\\(m+n-4)\mR^W\quad|\quad0\\-\delta_\phi^W\mR^C+mv^{-2}\nabla^W\tilde J\wedge\tilde J\end{pmatrix} ,
\end{equation}
as desired.
\end{proof}

To conclude this section, we observe that the Weyl tractor annihilates $W$-parallel sections of $\mT^W$.

\begin{prop}
\label{prop:WW_annihilates_parallel}
Let $(M^n,c,v^md\nu)$ be a SCMS with \charconstant\ $\mu$.  Suppose that $I\in\mT^W$ is such that $\nabla^WI=0$.  Then the contraction of $I$ into any component of $W^W$ vanishes.
\end{prop}

\begin{proof}

Since $\nabla^WI=0$, it follows immediately that $\mR^W(I)=0$ and $\nabla^W\mR^W(I)=0$.  Since $I\in\mT^W$, we have additionally that $\lp I,\tilde J\rp=0$ and $\lp I,\nabla^W\tilde J\rp=0$.  From the definition~\eqref{eqn:WC_coords} of $W^W$, we thus see that $W^W(\cdot,\cdot,I,\cdot)=0$.  By Corollary~\ref{cor:algebraic_tractor}, it follows that the contraction of $I$ into any component of $W^W$ vanishes, as desired.
\end{proof}
\section{Tractor Obstructions}
\label{sec:result}

We now turn to adapting the work of Gover and Nurowski~\cite{GoverNurowski2006} to our setting.  To that end, we will need two different formulations of what it means for a smooth conformal measure space to be generic.

\begin{defn}
A SCMS $(M^n,c,v^md\nu)$ with \charconstant\ $\mu$ is \emph{generic} if
\begin{enumerate}
\item $A^W\colon\Lambda^2T^\ast M\to\Lambda^2T^\ast M$ is injective,
\item $A^W\colon\End(TM)\to\End(TM)\oplus\Lambda^3T^\ast M\otimes TM$ given by
\[ T \mapsto \left( \lp A^W,T\rp, \sum_{i=1}^n A^W(e^i)\wedge T(e^i) \right) \]
is injective, where $\{e^i\}$ is an orthonormal basis for $TM$.
\end{enumerate}
\end{defn}

\begin{defn}
A SCMS $(M^n,c,v^md\nu)$ with \charconstant\ $\mu$ is \emph{weakly generic} if $A^W\colon TM\to\Lambda^2T^\ast M\otimes TM$ is injective.
\end{defn}

When the data of a SCMS and its \charconstant\ are clear from context, we shall merely say that $M$ is (weakly) generic.

\begin{remark}
These definitions are made \emph{pointwise}; i.e.\ we are viewing the maps $A^W$ as homomorphisms of vector bundles, not of their space of smooth sections (cf.\ Remark~\ref{rk:singular_obstruction}).
\end{remark}

In~\cite{GoverNurowski2006}, a SCMS is said to be \emph{$\Lambda^2$-generic} if $A^W\colon\Lambda^2T^\ast M\to\Lambda^2T^\ast M$ is injective.  This assumption was used by Listing~\cite{Listing2006} in his work on sharp metric obstructions to the existence of Einstein metrics in a conformal class, but will not be used here.  However, there is a nice chain of successively more general classes of SCMS,
\[ \{\text{generic SCMS}\} \subset \{\text{$\Lambda^2$-generic SCMS}\} \subset \{\text{weakly generic SCMS}\} . \]

Returning to the problem at hand, we observe that for generic SCMS, the weighted Weyl tractor gives a sharp obstruction to the existence of a quasi-Einstein metric.

\begin{thm}
\label{thm:obstruction_tractor_curvature}
Let $(M^n,c,v^md\nu)$ be a simply connected generic SCMS with \charconstant\ $\mu$.  There exists a positive quasi-Einstein scale $u\in\mE[1]$ if and only if there is a section $I\in\mT^W$ such that $\lp X,I\rp>0$, $W^W(I)=0$, and $\mR^W(I)=0$.
\end{thm}

\begin{remark}
If $m+n-4\not=0$, it is enough to assume $W^W(I)=0$.  Also, Theorem~\ref{thm:obstruction_tractor_curvature} is really a local result; the simply connected assumption is only used to conclude that a closed two-form is exact.  This is a feature of all of our results.
\end{remark}

\begin{proof}

By Theorem~\ref{thm:equivalence} and Proposition~\ref{prop:WW_annihilates_parallel}, all that remains to prove is that if there is a $I\in\mT^W$ such that $\lp X,I\rp>0$, $W^W(I)=0$, and $\mR^W(I)=0$, then there exists a $W$-parallel tractor $\hat I\in\mT^W$ with $\lp X,\hat I\rp>0$.  To that end, fix a scale $g\in c$ and denote $I=(\rho,\omega,\sigma)$.  Since $\sigma>0$ and $\mR^W(I)=0$, we have that
\begin{equation}
\label{eqn:k}
0 = dP^W - A^W(\sigma^{-1}\omega) \in \Lambda^2T^\ast M\otimes TM .
\end{equation}
Set $K=\sigma^{-1}\omega$.  First, taking the trace of~\eqref{eqn:k} and using~\eqref{eqn:trdP^C} and~\eqref{eqn:trA^C}, we have that
\begin{equation}
\label{eqn:trk}
0 = \tr dP^W + \tr A^W(K) = m(\sigma v)^{-1}\lp I,\nabla^W\tilde J\rp .
\end{equation}
Second, taking the weighted divergence in the $TM$ component of~\eqref{eqn:k} and using~\eqref{eqn:divA^C} and~\eqref{eqn:divdP^C}, it follows that
\begin{equation}
\label{eqn:divk}
0 = -m\sigma^{-1}v^{-2}\nabla v\wedge \lp I,\nabla^W\tilde J\rp - A^W(dK) + (m+n-3)dP^W(\cdot,\cdot,K) .
\end{equation}
Using~\eqref{eqn:k} again implies that the last summand vanishes, while~\eqref{eqn:trk} implies that the first summand vanishes.  Thus $A^W(dK)=0$, whence the genericity assumption implies that $K=df$ for some function $f\in\mE[0]$.  In particular, this implies that the tractor $\hat I=e^f\sigma^{-1}I$ is such that $\lp X,\nabla^W\hat I\rp=0$.  Moreover, working instead in the scale $e^{-2f}g,$ we may write $\hat I=(\rho,0,1)$, whereupon~\eqref{eqn:k} implies that $dP^W=0$.

Next, since $\mR^W$ is the curvature of $\nabla^W$, it automatically satisfies the Bianchi identity $d^W\mR^W=0$.  In particular, taking the exterior derivative of the equation $\mR^W(\hat I)=0$, it holds that
\[ 0 = \mR^W \wedge \nabla^W\hat I \in \Lambda^3T^\ast M\otimes \mT, \]
where here the wedge product includes a contraction in the tractor components in the obvious way.  Hence
\begin{equation}
\label{eqn:gen1}
0 = \lp Z, \mR^W\wedge\nabla^W\hat I\rp = \sum_{i=1}^n A^W(e^i)\wedge (P^W+\rho g)(e^i) \in \Lambda^2T^\ast M\otimes TM ,
\end{equation}
where $\{e^i\}$ is an orthonormal basis for $TM$.

Finally, $W^W(\hat I)=0$ also implies that
\[ (\delta_\phi^W\mR^W-mv^{-2}\nabla^W\tilde J\wedge\tilde J)(\hat I) = 0 . \]
Since $\hat I\in\mT^W$ and~\eqref{eqn:trk} holds, we have that $\delta_\phi^W\mR^W(\hat I)=0$.  It then follows from~\eqref{eqn:divphi_mR^C} and the fact that $dP^W=0$ that
\begin{equation}
\label{eqn:gen2}
0 = \lp Z, \delta^W\mR^W(\hat I)\rp = \lp A^W, P^W + \rho g\rp .
\end{equation}
Since $M$ is generic, \eqref{eqn:gen1} and~\eqref{eqn:gen2} together imply that $P^W+\rho g=0$ in the scale $e^{-2f}g$.  Thus $\nabla^W\hat I=\beta\otimes X$ for some one-form $\beta$.  By Lemma~\ref{lem:cute_algebraic_trick}, we thus have $\nabla^W\hat I=0$, as desired.
\end{proof}

If we are instead content to focus on the curvature of $\mR^W$ and its derivative $\nabla^W\mR^W$ --- which contains more information than $\delta_\phi^W\mR^W$ --- we can establish a sharp obstruction under the assumption that $M$ is weakly generic (cf.\ \cite[Theorem~3.4]{GoverNurowski2006}).

\begin{thm}
\label{thm:obstruction_tractor}
Let $(M^n,c,v^md\nu)$ be a simply connected weakly generic SCMS with \charconstant\ $\mu$.  There exists a positive quasi-Einstein scale $u\in\mE[1]$ if and only if there exists a nonvanishing tractor $I\in\mT^W$ such that
\begin{subequations}
\label{eqn:1-jet}
\begin{align}
\label{eqn:1-jet-0}\mR^W(I) & = 0 \\
\label{eqn:1-jet-1}\nabla^W\mR^W(I) & = 0 .
\end{align}
\end{subequations}
\end{thm}

\begin{proof}

Again by Theorem~\ref{thm:equivalence}, we only need to prove that the existence of a nonvanishing $I\in\mT^W$ satisfying~\eqref{eqn:1-jet} implies the existence of a $W$-parallel tractor $\hat I\in\mT^W$ satisfying $\lp X,\hat I\rp>0$.  To that end, denote $I=(\rho,\omega,\sigma)$, and suppose that $p\in M$ is such that $\sigma(p)=0$.  Then~\eqref{eqn:curvature_vector} and~\eqref{eqn:1-jet-0} imply that $A^W(\omega\rv_p)=0$, whence, by the genericity assumption, $\omega\rv_p=0$.  Since $I\not=0$, it holds that $\rho(p)\not=0$.  On the other hand, \eqref{eqn:1-jet-1} implies that at $p$,
\begin{equation}
\label{eqn:step1}
0 = \nabla^W\mR^W(\rho\otimes X) = -\rho(dP^W, A^W, 0) .
\end{equation}
In particular, $A^W\rv_p=0$, a contradiction.  Thus $\lp X,I\rp>0$.

Next, observe that since $M$ is weakly generic, the kernel of $\mR^W$ as a map $\mR^W\colon\mT\to\Lambda^2T^\ast M\otimes\mT$ has dimension at most two; this follows immediately from the fact that $\{\mR^W\cdot(0,e^i,0)\}\subset\Lambda^2T^\ast M\otimes\mT$ is linearly independent for any basis $\{e^i\}$ of $TM$.  However, by~\eqref{eqn:curvature_vector} and~\eqref{eqn:1-jet-0}, we know that $X$ and $I$ are in the kernel of $\mR^W$, and thus form a basis for $\ker\mR^W$.  On the other hand, \eqref{eqn:1-jet-0} and~\eqref{eqn:1-jet-1} imply that $\mR^W(\nabla^W I)=0$, whence
\[ \nabla^W I = \alpha\otimes I + \beta\otimes X \]
for some $\alpha,\beta\in T^\ast M$.  Applying $d^W$ and using~\eqref{eqn:1-jet-0} again, we find that $d\alpha=0$, whence $\alpha=df$ for some function $f$.

Consider now the tractor $\hat I := e^{-f}I$.  We still have that $\lp\hat I,\tilde J\rp=0$, while differentiating yields
\[ \nabla^W\hat I = e^{-f}\beta\otimes X . \]
Hence, by Lemma~\ref{lem:cute_algebraic_trick}, $\nabla^W\hat I=0$.
\end{proof}

\begin{remark}
\label{rk:singular_obstruction}
As observed by Gover~\cite{Gover2006}, we can weaken the genericity assumption by assuming only that $A^W\colon\Gamma(TM)\to\Gamma(\Lambda^2T^\ast M\otimes TM)$ is injective and still recover a sharp obstruction, provided we allow for scale singularities.  More precisely, under this assumption, one can still conclude the existence of a $W$-parallel tractor $I\in\mT^W$, but it is not necessarily the case that $\lp X,I\rp>0$.
\end{remark}

As an immediate corollary, we have the following sharp obstruction to the existence of quasi-Einstein metrics on weakly generic SCMS depending only on the $1$-jet of $\mR^W$, which moreover does not require finding a tractor satisfying certain conditions.

\begin{cor}
\label{cor:obstruction_tractor}
Let $(M^n,c,v^md\nu)$ be a simply connected weakly generic SCMS with \charconstant\ $\mu$.  Define
\[ E = \big(\Lambda^2T^\ast M\otimes\mT\big) \oplus \big(T^\ast M\otimes\Lambda^2T^\ast M\otimes\mT) \]
and consider the map $\Phi\colon\mT^W\to E$ given by
\begin{equation}
\label{eqn:linear_map}
\Phi(I) = \left( \mR^W(I), \nabla^W\mR^W(I) \right) .
\end{equation}
Then the induced map
\begin{equation}
\label{eqn:multilinear_map}
\Lambda^{n+1}\Phi\colon \Lambda^{n+1}\mT^W\to\Lambda^{n+1}E
\end{equation}
is the zero map if and only if there is a quasi-Einstein scale $u\in\mE[1]$.
\end{cor}

\begin{remark}
The condition that~\eqref{eqn:multilinear_map} be the zero map can easily be rephrased as the requirement that a system of natural conformal invariants depending only on $\mR^W$ and $\nabla^W\mR^W$, or equivalently on the two-jet of $A^W$, all vanish (cf.\ \cite{GoverNurowski2006}).
\end{remark}

\begin{proof}

By Theorem~\ref{thm:obstruction_tractor}, there exists a positive quasi-Einstein scale $u\in\mE[1]$ if and only if $\Phi$ is not injective.  But since $\dim\mT^W=n+1$, $\Phi$ is not injective if and only if $\Lambda^{n+1}\Phi$ is the zero map.
\end{proof}
\section{Sharp Tensorial Obstructions}
\label{sec:appl}

Using the results of the previous section as motivation, we can also find tensorial obstructions to the existence of quasi-Einstein metrics.  In particular, it is possible to determine candidates for the quasi-Einstein potential, should it exist.  In this way, we can find an invariant defined on weakly generic smooth conformal measure spaces which vanishes if and only if the space admits a quasi-Einstein metric, generalizing results of~\cite{BartnikTod2006,GoverNurowski2006,KNT1985}.  Moreover, these invariants can be extended to the interesting cases of static metrics and gradient Ricci solitons, which are not entirely covered by the results of the previous section.

\subsection{An invariant motivated by Section~\ref{sec:result}}
\label{sec:appl/tensor}

From Theorem~\ref{thm:obstruction_tractor}, we know that on a weakly generic SCMS $(M^n,c,v^md\nu)$ with \charconstant\ $\mu$, there are (local) sharp obstructions to the existence of a quasi-Einstein scale.  On the other hand, from the proof of Theorem~\ref{thm:obstruction_tractor_curvature}, we know that the key equation which gives the candidate $W$-parallel section of $\mT^W$ is~\eqref{eqn:k}, which, for convenience, we recall as
\begin{equation}
\label{eqn:k2}
0 = dP^W - A^W(K) .
\end{equation}
This is a purely tensorial expression, and in light of the proof of Theorem~\ref{thm:obstruction_tractor_curvature}, one might wonder if it is possible to solve for $K$ and produce a new tensorial expression which vanishes if and only if $M$ is quasi-Einstein.  In the special case of conformally Einstein metrics, this has already been carried out~\cite{GoverNurowski2006,KNT1985}, and as we shall see, their method can also be carried out in our setting.

The key observation is that since $M$ is weakly generic, it is possible to construct an ``inverse'' for $A^W$, and thereby solve for $K$.  Of course, since we are viewing $A^W$ as a homomorphism from $TM$ to $F:=\Lambda^2T^\ast M\otimes TM$, there are a number of ways to do this.  In order to simplify the exposition, we shall only discuss one possible method which works for Riemannian manifolds, and refer the reader to~\cite{GoverNurowski2006,KNT1985} for other possibilities.  To that end, we introduce the operator
\[ \vee \colon \Hom(TM,F)\times\Hom(TM,F) \to \End(TM) \]
by
\[ \lp(A\vee B)(x),y\rp = \lp A(x),B(y)\rp_F \]
for all $x,y\in TM$, and set $\check{A}=A\vee A$ --- using the metric to identify $\End(TM)\cong T^\ast M\otimes T^\ast M$, this is precisely the operator appearing in~\cite[1.131]{Besse}.  Since we have assumed that $g$ is Riemannian, it is clear that $\check{A}^W$ is invertible if and only if $M$ is weakly generic.  Thus we may define $D^W\in\Hom(TM,F)$ by
\[ D^W = A^W \circ (\check{A}^W)^{-1} . \]
Clearly $D^W\vee A^W=\id$, so $D^W$ is a desired inverse of $A^W$.  In particular,
\begin{equation}
\label{eqn:k3}
\lp dP^W - A^W(K), D^W(x)\rp_F = \lp dP^W, D^W(x)\rp_F - \lp K,x\rp .
\end{equation}

Now, since $M$ is weakly generic, it is clear that if $K_1,K_2\in TM$ are both such that~\eqref{eqn:k2} holds, then $K_1=K_2$.  Hence~\eqref{eqn:k3} implies that if~\eqref{eqn:k2} holds, then
\begin{equation}
\label{eqn:k4}
K = \lp dP^W, D^W(\cdot)\rp_F^\sharp,
\end{equation}
where $\sharp\colon T^\ast M\to TM$ is the usual ``musical'' isomorphism defined by $g$.  Of course, \eqref{eqn:k2} is not enough to conclude that $K$ comes from a quasi-Einstein potential.  However, the above discussion implies that $K$ comes from a positive quasi-Einstein scale with \charconstant\ $\mu$ if and only if the tensor
\begin{equation}
\label{eqn:G}
\begin{split}
G & = \Ric_\phi^m - \frac{R_\phi^m+m\mu v^{-2}}{m+n}g \\
& \quad + (m+n-2)\left(\nabla K + K\otimes K - \frac{1}{m+n}(\delta_\phi K+\lv K\rv^2)g\right)
\end{split}
\end{equation}
vanishes.  More precisely, if $u\in\mE[1]$ is a positive quasi-Einstein scale, then $G=0$ and $K=\nabla\log u$ by the uniqueness of solutions to~\eqref{eqn:k2}.  Conversely, if $G$ vanishes, then necessarily $\nabla K$ is symmetric.  In particular, locally we have $K=df$, and it follows from the vanishing of~\eqref{eqn:G} that $u=e^f\in\mE[1]$ is a quasi-Einstein scale with \charconstant\ $\mu$.

Finally, we recall that $A^W$ is the projecting part of $\mR^W$, and hence is conformally invariant.  By construction, this implies that $D^W$ is conformally invariant.  While $K$ is not conformally invariant, it satisfies a simple conformal transformation rule.  In particular, it is easy to check that~\eqref{eqn:G} is conformally invariant.  Thus, the tensor $G\in T^\ast M\otimes T^\ast M[0]$ is a well-defined tensor associated to a weakly generic SCMS $(M^n,g,v^md\nu)$ with \charconstant\ $\mu$.  Summarizing, we have established the following theorem.

\begin{thm}
\label{thm:obstruction}
Let $(M^n,c,v^md\nu)$ be a simply connected weakly generic SCMS with \charconstant\ $\mu$, and let $G$ be as in~\eqref{eqn:G}.  Then $G$ vanishes if and only if $M$ admits a positive quasi-Einstein scale $u\in\mE[1]$.
\end{thm}
\subsection{Obstructions for gradient Ricci solitons}
\label{sec:appl/grs}

As discussed in Section~\ref{sec:bg/qe}, an important aspect of the study of the conformal properties of a SMMS $(M^n,g,v^m\dvol)$ is that it makes sense when $\lv m\rv=\infty$, provided one writes $v^m=e^{-\phi}$ and rewrites the conformal factor in~\eqref{eqn:scms} as $u^{m+n-2}=e^f$ for $m$ finite and then takes the limit $\lv m\rv\to\infty$.  In this way, it is a trivial task to yield tensorial obstructions to the existence of gradient Ricci solitons from Theorem~\ref{thm:obstruction}, despite the fact that the tractorial obstructions do not make any sense.  The benefit to this approach is that it reveals immediately what are the correct analogues of the Weyl and the Cotton tensors if one were to try to generalize the tensorial perspectives of~\cite{BartnikTod2006,GoverNurowski2006,KNT1985}, which is not otherwise immediate obvious.

To make this precise, we need only understand the limiting behavior of the tensors $A^W$ and $dP^W$, as well as recall how the vector field $K$ is motivated.  Fix a constant $\mu\in\bR$ and a Riemannian manifold $(M^n,g)$, and regard them as the SMMS $(M^n,g,\dvol,\infty)$ with \charconstant\ $\mu$.  This is clearly the limit of the SMMS $(M^n,g,\dvol,m)$ with \charconstant\ $\mu$ as $m\to\infty$, and so one ask about the limiting behavior of the obstruction found in Theorem~\ref{thm:obstruction}.  Note that here we are fixing the scale $v=1$ in Theorem~\ref{thm:obstruction}, which is consistent with the usual problem of determining whether a given Riemannian manifold is a gradient Ricci soliton.

By way of motivation, suppose first that we can choose our sequence $m\to\infty$ such that for each $m$, the SMMS $(M^n,g,\dvol,m)$ admits a quasi-Einstein scale $u_m\in C^\infty(M)$ (depending on $m$) with \charconstant\ $\mu$.  Then the vector field $K$ defined in Theorem~\ref{thm:obstruction} is defined by
\[ K = \nabla\log u_m = \frac{1}{m+n-2}\nabla f_m, \]
where we define $u_m^{m+n-2}=e^{f_m}$, and we have that
\[ \Ric + \nabla^2 f_m - \frac{1}{m+n-2}df_m\otimes df_m - \frac{1}{m+n-2}\left(\Delta f_m - \lv\nabla f_m\rv^2\right)g = \mu g . \]
In particular, if we know additionally that the functions $f_m\to f\in C^\infty(M)$, then $(M^n,g,e^{-f}\dvol,\infty)$ is a quasi-Einstein SMMS with \charconstant\ $\mu$.

The above discussion tells us that, if we want to take the limit $m\to\infty$ in Theorem~\ref{thm:obstruction}, we should instead consider the vector field $\hat K=(m+n-2)K$ on the SMMS $(M^n,g,\dvol,m)$, so that~\eqref{eqn:k2} reads
\begin{equation}
\label{eqn:kinfty}
0 = (m+n-2)dP^W - A^W(\hat K) .
\end{equation}
Next, we observe that the tensors $P^W$, $A^W$, and $dP^W$ are such that
\begin{align}
\label{eqn:A^W_asymptotics} A^W & = \Rm + O\left(m^{-1}\right) \\
\label{eqn:dP^W_asymptotics} (m+n-2)dP^W & = d\Ric + O\left(m^{-1}\right) .
\end{align}
In particular, this tells us that the natural analogues of the Weyl tensor and the Cotton tensor in deriving obstructions in the spirit of~\cite{GoverNurowski2006,KNT1985} to the existence of quasi-Einstein metrics are the Riemann curvature tensor $\Rm$ and the derivative $d\Ric$ of the Ricci tensor, respectively.  This leads us to the following result.

\begin{thm}
\label{thm:grs_obstruction}
Let $(M^n,g)$ be a simply connected Riemannian manifold and suppose that the Riemann curvature $\Rm\in\Hom(TM,F)$ is injective, where $F$ is as in Section~\ref{sec:appl/tensor}.  Define $D\in\Hom(TM,F)$ by
\[ D = \Rm\circ\check{\Rm}^{-1} \]
and define the vector field $K$ by
\begin{equation}
\label{eqn:kinfty_defn}
K = \lp d\Ric, D(\cdot)\rp_F^\sharp .
\end{equation}
Then there exists a function $f$ such that $\Ric_f^\infty=\mu g$ if and only if
\begin{equation}
\label{eqn:infty_condn}
0 = \Ric + \nabla K - \mu g .
\end{equation}
\end{thm}

\begin{proof}

If~\eqref{eqn:infty_condn} holds, then necessarily $dK=0$.  Hence $K=df$ for some $f\in C^\infty(M)$, and so $\Ric_f^\infty=\mu g$.

Conversely, if $\Ric_f^\infty=\mu g$, then taking the exterior derivative of this equation yields
\[ 0 = d\Ric - \Rm(\nabla f) . \]
However, by construction~\eqref{eqn:kinfty_defn}, $K$ is a solution to
\[ 0 = d\Ric - \Rm(K) . \]
Since $\Rm$ is injective, this implies that $K=\nabla f$, and in particular, that~\eqref{eqn:infty_condn} holds.
\end{proof}
\subsection{Obstructions for static metrics}
\label{sec:appl/static}

An original motivation for the introduction of the tractor calculus to the study of quasi-Einstein metrics was to better understand certain similarities between static metrics and Poincar\'e-Einstein metrics.  In the context of this article, such a similarity can be seen between the obstructions to the existence of a static potential found by Bartnik and Tod~\cite{BartnikTod2006} and the obstruction to the existence of a conformally Einstein scale found by Kozameh, Newman and Tod~\cite{KNT1985} and later generalized by Gover and Nurowski~\cite{GoverNurowski2006}.  An important distinction here is that while the latter obstructions exist in all dimensions, the obstructions found by Bartnik and Tod only exist in dimension three, and it is not clear from their work how to extend their obstructions to higher dimensions.  The purpose of this section is to overcome this problem by once again using Theorem~\ref{thm:obstruction} to determine what are the right curvature tensors to use to find the higher dimensional analogue of the Bartnik-Tod obstruction tensor.

To clarify our terminology, by a \emph{static metric} on a smooth manifold $M^n$ we mean a metric $g$ such that $(M^n,g)$ can be realized as a time-symmetric spacelike hypersurface in a static $(n+1)$-dimensional spacetime with cosmological constant $\lambda$.  From an intrinsic point of view, this means that there exists a nonnegative function $v\in C^\infty(M)$, called the \emph{static potential}, such that
\begin{equation}
\label{eqn:std_static}
\begin{split}
0 & = v\Ric - \nabla^2 v + \Delta v\,g \\
R & = (n-1)\lambda .
\end{split}
\end{equation}
The corresponding static spacetime is the Lorentzian warped product $(\bR\times M, \og=-v^2dt^2\oplus g)$, which satisfies $\Ric(\og)=\lambda\og$.

In the language of the present article, if $g$ is a static metric on $M^n$ with static potential $v$, then the SMMS $(M^n,g,v^1\dvol)$ is a quasi-Einstein SMMS with quasi-Einstein constant $\lambda$ and \charconstant\ zero.  Since $m=1$, such a SMMS can be treated using tractor methods.  In particular, Theorem~\ref{thm:obstruction} can be applied as follows: Given a simply connected weakly generic SMMS $(M^n,g,v^1\dvol)$ with \charconstant\ zero, Theorem~\ref{thm:obstruction} yields a symmetric $(0,2)$-tensor which vanishes if and only if $(M^n,g,v^1\dvol)$ is conformally static, making no restriction on $\lambda$.

From the classical point of view in general relativity, this perspective is somewhat lacking, as it does not allow one to determine whether a given metric $g$ is itself a static metric.  This can be overcome by Corollary~\ref{cor:duality}, which implies that if $(M^n,g)$ is a static manifold with static potential $v$ and cosmological constant $\lambda$, then the SMMS $(M^n,g,1^{1-n}\dvol_g)$ is such that $v\in C^\infty(M)$ is a quasi-Einstein scale with quasi-Einstein constant zero and \charconstant\ $\lambda$.  In other words, addressing the question of whether a given Riemannian manifold $(M^n,g)$ admits a static potential is equivalent to determining whether the SMMS $(M^n,g,1^{1-n}\dvol_g)$ with \charconstant\ $\lambda$ admits a quasi-Einstein scale.

Unfortunately, one \emph{cannot} use Theorem~\ref{thm:obstruction} to determine whether the SMMS $(M^n,g,1^{1-n}\dvol_g)$ with \charconstant\ $\lambda$ admits a quasi-Einstein scale, because the dimensional constant $1-n$ is one of the disallowed values.  The reason this value is disallowed is because the weighted Schouten tensor of the SMMS $(M^n,g,1^{2-m-n}\dvol_m)$ with \charconstant\ $\lambda$ is
\begin{equation}
\label{eqn:weighted_schouten_prestatic}
P^W = -\frac{1}{m}\left(\Ric + \frac{R-(m+n-1)\lambda}{2(m-1)}g\right),
\end{equation}
which is singular for $m=1$.  Nevertheless, by imposing the additional (physically reasonable) assumption that the quasi-Einstein scale $v$ have quasi-Einstein constant zero, it will be possible to overcome this difficulty if one is willing to work only in a tensorial language.

To motivate our result, let us first give a heuristic overview of how to overcome the singularity at $m=1$ in~\eqref{eqn:weighted_schouten_prestatic} using the additional requirement that any quasi-Einstein scale have quasi-Einstein constant zero.  First, fix a SMMS $(M^n,g,1^{2-m-n}\dvol_g)$ with \charconstant\ $\lambda$.  If $v$ is such a quasi-Einstein scale with quasi-Einstein constant $\mu$, then
\begin{equation}
\label{eqn:norm_v_prestatic}
\frac{\mu}{m-1} = \left|\frac{1}{m-2}\bD^W v\right|^2 = \frac{1}{m-2}\left(v\Delta v-\frac{R-(m+n-2)\lambda}{2(m-1)}v^2\right) + \lv\nabla v\rv^2 .
\end{equation}
Multiplying through by $m-1$, we see that this still makes sense when $m=1$, where it forces $\mu=0$ and $R=(n-1)\lambda$ if $v$ is not constant; in other words, the constraint $R=(n-1)\lambda$ in~\eqref{eqn:std_static} is the condition which forces $\mu=0$.  In this way, one is tempted to na\"ively \emph{declare} the weighted Schouten tensor of $(M^n,g,1^{1-n}\dvol_g)$ with \charconstant\ $\lambda$ to be $P^W=-\Ric$, which arises by ignoring the last term in~\eqref{eqn:weighted_schouten_prestatic} via this cancellation.  Likewise, one also wants to \emph{declare} $A^W=\Rm+\Ric\wedge g$.  Of course, one cannot actually make sense of this via limits --- any such statement would require including the quasi-Einstein scale $v$ in the definition of $P^W$, which is undesirable because $v$ is not assumed \emph{a priori} to exist --- but it nevertheless leads one effortlessly to the following generalization of the obstruction found by Bartnik and Tod~\cite{BartnikTod2006}.

\begin{thm}
\label{thm:static_obstruction}
Let $(M^n,g)$ be a simply connected Riemannian manifold.  Define the tensor $B\in\Hom(TM,F)$ by
\[ B(x) = A^W(x) - \frac{1}{n-1}\tr A^W(x) \wedge g \]
for all $x\in TM$, where $A^W=\Rm+\Ric\wedge g$, and suppose that $B$ is injective.  Define the tensor $D\in\Hom(TM,F)$ by
\[ D = B\circ \check{B}^{-1}, \]
and set
\begin{equation}
\label{eqn:static_k}
K=-\lp d\Ric,D(\cdot)\rp_F^\sharp.
\end{equation}
The $(M^n,g)$ admits a static potential with cosmological constant $\lambda$ if and only if the scalar curvature $R=(n-1)\lambda$ and the $(0,2)$-tensor
\begin{equation}
\label{eqn:static_G}
G = \Ric - \nabla K - K\otimes K - \lambda g
\end{equation}
vanishes.
\end{thm}

\begin{proof}

First suppose that $R=(n-1)\lambda$ and~\eqref{eqn:static_G} vanishes.  Then $dK=0$, whence there is a function $f\in C^\infty(M)$ such that $K=df$.  Set $v=e^{f}$, so that~\eqref{eqn:static_G} implies that
\[ 0 = v\Ric - \nabla^2v - \lambda v g . \]
Taking the trace and using the assumption on $R$ implies that $\Delta v+\lambda v=0$, yielding~\eqref{eqn:std_static}.

Conversely, suppose that~\eqref{eqn:std_static} holds.  Differentiating, we see that
\begin{equation}
\label{eqn:static_c}
\begin{split}
0 & = d\Ric + v^{-1}(\Rm+\Ric\wedge g)(\nabla v) - v^{-1}(\Ric+\lambda g)(\nabla v)\wedge g \\
& = -\,dP^W + A^W(\nabla\log v) - \frac{1}{n-1}\tr A^W(\nabla\log v)\wedge g,
\end{split}
\end{equation}
where in the second line we have written $P^W=-\Ric$ and $A^W=\Rm-P^W\wedge g$ as suggested by the heuristic above.  However, by construction we also have that
\[ 0 = -dP^W + A^W(K) - \frac{1}{n-1}A^W(K)\wedge g, \]
whence, by the assumption that $B$ is injective, we have $K=\nabla\log v$.  Taking the trace of~\eqref{eqn:std_static} implies that $\Delta v+\lambda v=0$, whence follows the vanishing of~\eqref{eqn:G}.
\end{proof}

\begin{remark}

One could also define~\eqref{eqn:static_G} as the direct limit of~\eqref{eqn:G}, namely
\[ G = \Ric - \nabla K - K\otimes K + \left(\delta K + \lv K\rv^2\right)g, \]
which is also more directly related to the formulation~\eqref{eqn:std_static}.  We have opted instead to use the definition~\eqref{eqn:static_G} to parallel the obstruction found in~\cite{BartnikTod2006}.
\end{remark}

Note in particular that in Theorem~\ref{thm:static_obstruction} we do not use the curvature tensor $A^W$ to construct $G$, but rather some ``trace'' modification of it.  Unfortunately, we do not have a nice explanation for why this tensor $B$ should arise instead of $A^W$, other than the direct computation~\eqref{eqn:static_c}.  This is no doubt related to our heuristic motivation for the definition of $P^W$ and $A^W$, and especially the fact that the ``limits'' are not truly well-defined.

Let us conclude this section by showing that Theorem~\ref{thm:static_obstruction} gives exactly the obstruction found by Bartnik and Tod~\cite{BartnikTod2006} in dimension $n=3$, with the trivial generalization of allowing for a nonzero cosmological constant.  First, in dimension three, $B(K)$ simplifies to
\[ B(K) = \Ric_0(K)\wedge g + 2K\wedge\Ric_0 , \]
where $\Ric_0=\Ric-\frac{R}{3}g$ is the traceless Ricci curvature.  Letting $\{e^1,e^2,e^3\}$ be an orthonormal basis of eigenvectors of $\Ric_0$ with eigenvalues $\{\lambda_1,\lambda_2,\lambda_3\}$, we thus see that
\begin{equation}
\label{eqn:B3}
\begin{split}
B(K) & = (\lambda_2-\lambda_3)\lp K,e^1\rp \big(e^1\wedge e^2\otimes e^2+e^3\wedge e^1\otimes e^3\big) \\
& \quad + (\lambda_3-\lambda_1)\lp K,e^2\rp \big(e^1\wedge e^2\otimes e^1+e^2\wedge e^3\otimes e^3\big) \\
& \quad + (\lambda_1-\lambda_2)\lp K,e^3\rp \big(e^3\wedge e^1\otimes e^1+e^2\wedge e^3\otimes e^2\big) ,
\end{split}
\end{equation}
where we have used the fact $\tr\Ric_0=0$.  In particular, $B$ is injective if and only if the eigenvalues of $\Ric_0$, and hence $\Ric$, are distinct.  This is precisely the genericity assumption used by Bartnik and Tod~\cite{BartnikTod2006}.

\begin{remark}
By using the Hodge star $\ast\colon\Lambda^2T^\ast M\to\Lambda^1T^\ast M$ on $M^3$, \eqref{eqn:B3} gives exactly the coordinate expression relating the static potential to the Ricci and Cotton-York tensors which appears in~\cite{BartnikTod2006,Tod2000}.  From this, is it also easy to give an explicit formula for the gradient of the static potential in the Ricci eigenframe, and, with a little more work, an invariant expression for $K$ as in~\eqref{eqn:static_k} in dimension three; for details, see~\cite{BartnikTod2006}.
\end{remark}

\begin{remark}
As in~\cite{BartnikTod2006}, Theorem~\ref{thm:static_obstruction} can be regarded as yielding an algorithm for finding the static potential on a weakly generic Riemannian manifold, if it exists.  This is done by asking the following sequence of questions.
\begin{enumerate}
\item Is the scalar curvature constant?  If so, set $R=(n-1)\lambda$.
\item Is the vector field $K$ defined by~\eqref{eqn:static_k} locally a gradient?  If so, set $K=\nabla\log v$.
\item Does $v\Ric-\nabla^2v+\Delta v\,g=0$?  If so, $v$ is the static potential for $M$.
\end{enumerate}
\end{remark}
\section{A New Perspective on Hamilton's Matrix Harnack Inequality}
\label{sec:conclusion}

We conclude by observing a relationship between the weighted Weyl tractor $W^W$ and Hamilton's matrix Harnack inequality~\cite{Hamilton1993}.  To that end, first recall that Hamilton showed that if $(M,g(t))$ is a solution to the Ricci flow for which $g_0$ has nonnegative curvature operator, then
\begin{equation}
\label{eqn:harnack}
B(\alpha,\alpha) - 2\lp d\Ric, \Phi\otimes\alpha\rp + \Rm(\Phi,\Phi)\geq 0
\end{equation}
for all $g(t)$ and all $\alpha\in T^\ast M$, $\Phi\in\Lambda^2T^\ast M$, where
\begin{equation}
\label{eqn:hamilton_B}
B = \Delta\Ric - \frac{1}{2}\nabla^2 R + 2\lp\Rm,\Ric\rp - \Ric^2 + \frac{1}{2}\Ric .
\end{equation}
(cf.\ \cite{ChowChu1995}).  If one assumes that $(M,g(0))$ is a Ricci soliton, so that $g(t)$ is obtained from $g(0)$ by time-dependent dilations and diffeomorphisms only, then~\eqref{eqn:harnack} holds along $g(t)$ if and only if it holds for $g(0)$.  Hamilton raised the question of whether~\eqref{eqn:harnack} can be regarded as the statement that some bundle associated to the Ricci flow has nonnegative curvature, which has now been verified in a number of different ways~\cite{Cabezas-RivasTopping2009,ChowChu1995,Perelman1}.  Here, we demonstrate that the nonnegativity of~\eqref{eqn:harnack} for the metric $g=g(0)$ can be regarded as the nonnegativity of the weighted Weyl tractor associated to the SMMS $(M^n,g,\dvol,\infty)$ in the following sense:

\begin{thm}
\label{thm:harnack}
Let $(M^n,g)$ be a Riemannian manifold and suppose that there is a sequence $\{m_i\}\subset\bR$ with $m_i\to\infty$ such that the weighted Weyl tractors $W_i^W$ of the SMMS $(M^n,g,1^{m_i}\dvol)$ with \charconstant\ $\mu=-\frac{1}{2}$ are nonnegative.  Then~\eqref{eqn:harnack} holds for all $\alpha\in T^\ast M$, $\Phi\in\Lambda^2T^\ast M$.
\end{thm}

\begin{proof}

For convenience, fix $i$ and denote $m=m_i$.  Fix $\alpha\in T^\ast M$ and $\Phi\in\Lambda^2T^\ast M$.  Given any $\phi\in C^\infty(M)$ and $\beta\in T^\ast M$, form the tractor
\[ T = \begin{pmatrix} (m+n-4)\alpha\\\Phi\quad|\quad\phi\\\beta\end{pmatrix} \in \mA, \]
so that~\eqref{eqn:w_eval} yields
\begin{equation}
\label{eqn:w_nonneg}
\begin{split}
\frac{1}{m+n-4}W^W(T,T) & = (m+n-4)B^W(\alpha,\alpha) \\
& \quad - 2(m+n-4)\lp dP^W,\Phi\otimes\alpha\rp + A^W(\Phi,\Phi) .
\end{split}
\end{equation}
In addition to~\eqref{eqn:A^W_asymptotics} and~\eqref{eqn:dP^W_asymptotics}, we have that
\begin{equation}
\label{eqn:B^W_asymptotics}
(m+n-2)B^W = \delta d\Ric + \lp\Rm,\Ric\rp - \mu\Ric + O\left(m^{-1}\right) .
\end{equation}
Using the Weitzenb\"ock formula
\[ (d\delta+\delta d)T = \Delta T - T\circ\Ric + \lp\Rm,T\rp \]
for symmetric $(0,2)$-tensors together with the Bianchi identity $\delta\Ric=\frac{1}{2}dR$, we see that~\eqref{eqn:B^W_asymptotics} is equivalent to
\begin{equation}
\label{eqn:B2}
(m+n-2)B^W = \Delta\Ric - \frac{1}{2}\nabla^2 R + 2\lp\Rm,\Ric\rp - \Ric^2 - \mu\Ric + O\left(m^{-1}\right) .
\end{equation}
In particular, we may rewrite~\eqref{eqn:w_nonneg} as
\[ W^W(T,T) = B(\alpha,\alpha) - 2\lp d\Ric,\Phi\otimes\alpha\rp + \Rm(\Phi,\Phi) + O\left(m^{-1}\right) \]
for $B$ as in~\eqref{eqn:hamilton_B}, where we have now used that $\mu=-\frac{1}{2}$.  Since $W^W$ is nonnegative for all $i$, we may then take the limit $i\to\infty$ to yield the desired conclusion.
\end{proof}

In the proof, we have seen that the choice of \charconstant\ $\mu=-\frac{1}{2}$ only appeared to get the exact form~\eqref{eqn:hamilton_B}.  Since the choice $\mu=-\frac{1}{2}$ is analogous to searching for expanding gradient Ricci solitons in the sense of Section~\ref{sec:appl/grs}, this perspective fits in nicely with the perspectives of~\cite{Cabezas-RivasTopping2009,ChowChu1995,Hamilton1993}, where one always considered as a model an expanding gradient Ricci soliton in order to arrive at~\eqref{eqn:harnack}.

From the perspective of addressing Hamilton's original question~\cite{Hamilton1993}, Theorem~\ref{thm:harnack} is lacking in that it does not provide any hint as to how the Ricci flow influences this picture, nor does it suggest how one might actually prove Hamilton's result (cf.\ \cite{Cabezas-RivasTopping2009,ChowChu1995,Hamilton1993}).  It would be interesting to know if our perspective can be refined to meet these challenges.

\bibliographystyle{abbrv}
\bibliography{../bib}
\end{document}